\documentclass[nojournal,nonblindrev]{informs3} 

\DoubleSpacedXI 

\DeclareMathOperator*{\maximize}{maximize: }
\newcommand{\ttheta}{ \tilde{\theta} }

\newcommand{\bx}{ \mathbf{x} }
\newcommand{\by}{ \mathbf{y} }
\newcommand{\bz}{ \mathbf{z} }
\newcommand{\bxh}{ \mathbf{\hat{x}} }
\newcommand{\byh}{ \mathbf{\hat{y}} }

\usepackage{multirow}
\usepackage[table,xcdraw]{xcolor}
\usepackage{makecell}
\newcolumntype{L}[1]{>{\raggedright\let\newline\\\arraybackslash\hspace{0pt}}m{#1}}
\newcolumntype{C}[1]{>{\centering\let\newline\\\arraybackslash\hspace{0pt}}m{#1}}
\newcolumntype{R}[1]{>{\raggedleft\let\newline\\\arraybackslash\hspace{0pt}}m{#1}}
\usepackage{lipsum} 
\usepackage{tikz}

\usetikzlibrary{arrows.meta} 

\usepackage{endnotes}
\let\footnote=\endnote
\let\enotesize=\normalsize
\def\notesname{Endnotes}%
\def\makeenmark{$^{\theenmark}$}
\def\enoteformat{\rightskip0pt\leftskip0pt\parindent=1.75em
  \leavevmode\llap{\theenmark.\enskip}}

\usepackage{natbib}
\def\bibfont{\small}%
\usepackage{enumitem}
\usepackage{graphicx}      
\usepackage{caption}       
\usepackage{subfig}        
\usepackage{float}

\TheoremsNumberedThrough     
\ECRepeatTheorems

\EquationsNumberedThrough    

\usepackage{todonotes}

\begin{document}


\RUNAUTHOR{Stratman, Boutilier, and Albert}

\RUNTITLE{Ambulance Allocation for Patient-Centered Care}

\TITLE{Ambulance Allocation for Patient-Centered Care}

\ARTICLEAUTHORS{%
\AUTHOR{Eric G. Stratman}
\AFF{Department of Industrial and Systems Engineering,
University of Wisconsin-Madison, \EMAIL{egstratman@wisc.edu}}
\AUTHOR{Justin J. Boutilier}
\AFF{Telfer School of Management,
University of Ottawa, \EMAIL{boutilier@telfer.uOttawa.ca}}
\AUTHOR{Laura A. Albert}
\AFF{Department of Industrial and Systems Engineering,
University of Wisconsin-Madison, \EMAIL{laura@engr.wisc.edu}}
}

\ABSTRACT{
Problem definition: Emergency Medical Services (EMS) in the United States and similar systems typically utilize a single treatment pathway, transporting all patients to emergency departments (EDs), regardless of their actual care needs or preferences. Recent policy reforms have sought to introduce alternative treatment pathways to divert lower acuity patients from the ED, but operationalizing these options has proven difficult. This paper proposes a patient-centered EMS (PC-EMS) ambulance allocation model that supports multiple care pathways by aligning EMS responses with individual patient needs. Methodology/results: We develop a two-stage mixed-integer optimization framework that incorporates multiple dispatch and secondary assignment strategies which enable dynamic resource deployment. The model maximizes appropriate ED diversions while maintaining ambulance availability using a queueing-based availability constraint. We leverage national EMS data and machine learning to estimate dispatcher accuracy and diversion potential. Simulations across diverse geographic regions suggest that agencies can achieve up to 80\% of possible ED diversions by equipping only 15–25\% of their fleet with diversion-capable units. Adaptive dispatch strategies improve diversion rates by 3.4 to 8.6 times compared to conventional single-unit dispatch. Managerial implications: These results provide actionable guidance for PC-EMS implementation by quantifying the trade-off between equipment investment and operational coordination. Using the allocation model, agencies can strategically choose between upgrading fewer units with advanced dispatching protocols versus larger fleet investments with simpler operations. This evidence-based approach offers flexible pathways suited to different organizational capabilities and implementation readiness.
}%

\KEYWORDS{Facility location, emergency response, emergency medical services, machine-learning, queuing.} 

\maketitle
\section{Introduction}
Emergency Medical Services (EMS) in the United States, Canada, Australia, and other countries with rapid transport systems have traditionally followed a single-pathway model in which nearly all patients experiencing medical emergencies are transported to hospital emergency departments (EDs), regardless of their actual care needs or preferences. This one-size-fits-all approach, while operationally straightforward, not only contributes to ED overcrowding and rising healthcare costs but also neglects patient choice in their care journey. A growing body of research shows that many EMS patients could be appropriately treated in lower-acuity settings, often with better outcomes and at lower cost \citep{billings_emergency_2000}.

In response, policymakers have introduced initiatives that allow EMS providers to tailor transport decisions to a patient’s condition (i.e., their clinical needs or preferences). For example, the Emergency Triage, Treat, and Transport (ET3) model, launched by the Centers for Medicare and Medicaid Services, enabled EMS providers to treat patients in place (TIP) via telemedicine or transport them to alternative destinations (AD), such as urgent care centers. Although ET3 aligned with the goals of improving efficiency and reducing unnecessary ED visits, it was ultimately discontinued due to limited use of the alternative pathways \citep{centers_for_medicare__medicaid_services_emergency_2022}. These challenges suggest that transitioning EMS systems from a one-size-fits-all model to one that supports multiple treatment pathways requires more than just policy change. To prioritize patient-choice through ED diversions, EMS agencies need practical guidance on how to implement these changes effectively.

This paper addresses this need by presenting a patient-centered EMS (PC-EMS) allocation model designed to navigate the complexities of multiple treatment pathways. The model optimizes how ambulance units should be allocated to stations and assigned to patients, maximizing the number of eligible patients safely diverted from the ED. The approach provides a structured framework and insights for EMS agencies seeking to adopt flexible care pathways while maintaining ambulance availability and resource efficiency. 

\subsection{Challenges of Patient-Centered Emergency Medical Services}
The transition from traditional EMS models to patient-centered approaches requires addressing fundamental operational complexities that existing EMS allocation frameworks do not capture. While EMS allocation has been extensively studied, most models assume a single treatment pathway where all patients are transported to emergency departments, regardless of their specific needs. This framework fails to address three critical challenges that emerge when implementing PC-EMS systems with multiple treatment options.

First, not all EMS units are equipped to provide alternative treatment pathways. As highlighted in the ET3 program, capabilities such as telemedicine support or additional clinical training are necessary for TIP or AD referrals, and these capabilities are not universally available across all ambulances. This creates a meaningful distinction between traditional ambulances, which are limited to ED transport, and diversion-capable units, which have the necessary resources to support patient-centered pathways. While some prior EMS models account for heterogeneity in vehicle types (e.g., basic vs. advanced life support), they do not typically consider how those vehicle capabilities influence the treatment pathways available to the patient \citep{mclay_maximum_2009, boujemaa_multi-period_2020}.

Second, PC-EMS introduces a new decision point into the EMS response process. Initial dispatch decisions must be made before the patient’s condition is fully understood, based only on information attained during dispatcher screening. The actual treatment pathway is only revealed after on-scene assessment. This introduces uncertainty into initial vehicle assignments and may necessitate additional resources if the first responding unit lacks the capabilities to deliver the appropriate care. Although a few studies explore how EMS systems can adapt to updated information post-assessment \citep{marianov_hierarchical_2001, yoon_dynamic_2020}, this remains an underexplored area in the literature.

Third, the introduction of multiple treatment pathways leads to greater variability in service times. Calls involving TIP or AD transport may require more (or less) time than traditional ED transport due to the nature of the treatment or the travel times. Moreover, when vehicle substitution is required—i.e., dispatching a second unit because the initial one is not diversion-capable—service times can increase even further. Traditional EMS availability models often rely on fixed or homogeneous service time assumptions, which fail to account for these endogenous substitutions \citep{marianov_locationallocation_2002, boutilier_drone_2022}. In contrast, we extend existing constraints to accommodate this new source of variability, allowing us to better model resource availability in PC-EMS systems.

To address these operational challenges, we develop an approach for PC-EMS resource allocation, as outlined in the following section.

\subsection{PC-EMS Allocation Model}
We propose a PC-EMS resource allocation model to maximize the number of patients appropriately diverted from EDs while maintaining ambulance availability. To address the uncertainty in patient eligibility for alternative pathways and decision point after the on-scene patient assessment, we develop a two-stage mixed-integer optimization model.

In the first stage, the model determines the optimal allocation of EMS vehicles to stations and assigns these resources to patients based on the limited information available from dispatcher screening. This stage establishes the initial response configuration before patient conditions are fully understood. In the second stage, after EMS personnel conduct an in-person evaluation and determine the appropriate treatment pathway, the model updates resource assignments to align with the patient's now known eligible treatment options.

To effectively align vehicle capabilities with patient needs and maximize appropriate diversions, we explicitly incorporate two EMS assignment strategies that are commonly used in practice but have not been fully modeled within traditional EMS allocation frameworks: multiple response and secondary assignment. Multiple response involves dispatching more than one unit simultaneously when there is uncertainty about the required level of care, ensuring that at least one responding unit has the expertise necessary to assess the patient and determine the appropriate treatment pathway. Secondary assignment, implemented as a recourse mechanism in our model, enables additional resources to be deployed if the initial responding unit determines that different expertise is needed after patient assessment. These strategies create a flexible system capable of adapting to evolving information throughout the emergency response process.

Lastly, to prevent resource congestion from these strategies, we introduce a availability constraint based on an $M/G/d/d$ queuing framework. This constraint ensures that ambulance availability is maintained while balancing the trade-offs associated with deploying additional resources.

\subsection{Our Contributions}
The key contributions of our research include:
\begin{enumerate}
    \item \textbf{Two-Stage PC-EMS Allocation Model:} We develop the first mixed-integer optimization model that integrates multiple treatment pathways and dynamic resource deployment for PC-EMS systems. Unlike previous EMS allocation models that assume fixed treatment pathways, our two-stage approach explicitly captures how dispatch decisions evolve as patient information is revealed, incorporating both multiple response and secondary assignment strategies while using an $M/G/d/d$ queue framework to maintain system availability.
    \item \textbf{Machine Learning Framework for Uncertainty Estimation:} We develop a data-driven approach that quantifies two key uncertainties in PC-EMS: the likelihood of patients being eligible for non-ED pathways and the accuracy of dispatcher predictions. This framework enables realistic parameterization of the allocation model, allowing decision-makers to understand how pre-arrival information quality affects system performance.
    \item \textbf{System Simulation and Policy Insights:}  We evaluate our approach through simulation across diverse geographic environments. Our findings suggest that agencies can achieve 80\% of potential diversions with only 15-25\% of their fleet equipped with diversion capabilities. Multiple response and recourse dispatching strategies significantly improve outcomes when resources are limited, with agencies experiencing 3.4-8.6 times more diversions compared to single-unit dispatch. We also quantify the impact of dispatcher screening accuracy, providing practical guidance for PC-EMS implementation.
\end{enumerate}
The remainder of this paper is structured as follows. Section 2 reviews relevant literature on EMS allocation and predictive decision-making in emergency response. Section 3 presents our two-stage mixed-integer optimization model for PC-EMS allocation. Section 4 introduces our machine learning framework for estimating patient eligibility and dispatcher accuracy. Section 5 evaluates system performance through simulation using real-world EMS data. Section 6 discusses managerial and policy insights, and Section 7 concludes with key findings and directions for future research.
\section{Literature Review} \label{sec_lit_rev}
EMS system planning represents a rich body of research, and we refer readers to comprehensive reviews for broader field overviews \citep{stratman_uncertainty_2023}. Our work aligns most closely with EMS system research examining how diverse resources address heterogeneous patient populations.

In this domain, current approaches to emergency medical services recognize that different patient types present varying needs within ambulance systems. Early strategies position resources to ensure coverage by multiple unit types, acknowledging the different capabilities each resource offers \citep{marianov_probabilistic_1992, mandell_covering_1998}. Contemporary methods investigate how specific resources align with different patient requirements \citep{mclay_maximum_2009, boujemaa_multi-period_2020, nelas_locating_2021, yoon_dynamic_2021}. While these allocation models recognize distinct patient dispatching needs, few studies address how information evolves during emergency response. \cite{stratman_patient_2022} demonstrate how provider impressions develop as additional information becomes available throughout a response.

Our approach of assigning ambulances based on initial patient condition assessment and then adapting response when new information emerges is similar to \cite{ansari_approximate_2017} and \cite{yoon_dynamic_2020}, who examine multiple response strategies for managing patient need uncertainty, and \cite{marianov_hierarchical_2001}, who investigate resource actions similar to our secondary assignment mechanism. We extend these studies by integrating both approaches within a unified optimization framework that maximizes ED diversions while maintaining system availability.

Resource allocation under uncertainty constitutes another critical research area. Our study focuses particularly on uncertainty in ambulance availability. Traditional approaches incorporate availability measures directly into objective functions \citep{daskin_maximum_1983, batta_maximal_1989, restrepo_erlang_2009}. Alternative methods employ chance constraints to enforce probabilistic guarantees on system performance. \cite{marianov_probabilistic_1998} introduce frameworks ensuring server availability with specified levels, while \cite{boutilier_drone_2022} expand these methods to more general settings. Our work extends these approaches by considering availability with recourse actions through secondary assignment.

Regarding dispatcher accuracy, researchers assess dispatcher effectiveness in classifying patient urgency \citep{shah_derivation_2003, bohm_accuracy_2018}, while others explore machine learning approaches for improving non-critical patient need prediction \citep{eastwood_patient_2018, kang_artificial_2020}. Several operational studies examine the system-level impact of screening accuracy \citep{argon_priority_2009, mclay_model_2013}. Our research distinctively integrates screening accuracy metrics into strategic allocation, providing actionable insights on how dispatcher performance influences resource utilization and ED diversion rates across various dispatching strategies.

\section{Optimization Approach} \label{opt_mod}
We now present a two-stage allocation model with recourse to determine the optimal allocation strategy of EMS resources in a PC-EMS system. The model seeks to maximize patient diversions from the ED while maintaining adequate system-wide ambulance availability and other constraints on the EMS system.

In the first stage, EMS resources are allocated to stations, and an initial assignment strategy is defined, specifying the units that initially respond to a patient based on their location and the preliminary condition assessment from phone screening. In the second stage, after EMS personnel conduct an in-person evaluation and learn the patient’s true condition and appropriate treatment pathway, the model determines how the initially assigned units are utilized and whether to initiate secondary assignments, sending additional units in cases where the initial response is insufficient. 

\subsection{First-Stage Model} \label{sec_first_stage}
Consider an EMS system with several types of EMS units represented by the set $K$. We wish to locate $N_k$ of each type, $k \in K$, to a set of locations, $I$, to serve patients at a set of demand nodes, $J$. We let $\Theta$ be the set of patient conditions assessed from phone screening and note each patient has a believed condition, $\theta \in \Theta$, that helps the EMS dispatcher predict the treatment option the patient can utilize. 

Let the binary decision variables $z_{ikn} = 1$ when there are at least $n \in \{1,...,N_k\}$ EMS units of type $k \in K$ located at station $i \in I$ and $0$ otherwise. The binary decision variables $y_{ijk\theta} = 1$ indicate that a unit of type $k \in K$ located at station $i \in I$ is assigned to a patient at $j \in J$ with believed condition $\theta \in \Theta$. This assignment is made before the treatment pathways appropriate for a patient's actual condition is known. We impose a set of constraints on dispatching decisions, represented by $\mathcal{C}$. These constraints enforce commonly used performance criteria in ambulance systems, such as response-time threshold coverage, or average response time, or equity constraints. We formulate the first-stage model as follows:
\begin{align}
     \maximize_{\mathbf{y},\mathbf{z}}  \text{ } & D(\mathbf{y},\mathbf{z})
    \label{eq:objective} \\
    \text{subject to: } 
        & \sum_{i \in I}\sum_{n = 1}^{N_k } z_{ikn} \le N_k, 
    \quad \forall k \in K, \label{eq:constraint1}\\
      & z_{ikn} \le z_{ik(n-1)}, 
    \quad \forall i \in I, k \in K, n \in \{2,...,N_k\}, \label{eq:constraint2}\\
    & y_{ijk\theta } \le z_{ik1}, 
\quad \forall i \in I, j \in J, k \in K, \theta \in \Theta,  \label{eq:constraint3}\\  
    & \sum_{i \in I} \sum_{k \in K} y_{ijk\theta } \ge 1, 
\quad \forall j \in J, \theta \in \Theta,  \label{eq:constraint33}\\  
   &\mathbf{y} \in \mathcal{C},\label{eq:constraint4}\\
     &y_{ijk\theta } \in \{0,1\}, 
    \quad \forall i \in I, j \in J, k \in K, \theta \in \Theta, \label{eq:constraint5} \\
     &  z_{ikn}  \in \{0,1\},  \quad \forall i \in I, k \in K,n \in \{1,...,N_k\}. \label{eq:constraint6}
\end{align}

The objective function \eqref{eq:objective} maximizes the function, $D(\mathbf{y},\mathbf{z})$, which we define as the maximum number of patients diverted from the ED for given dispatching decisions, $\mathbf{y}$, and resource allocation decisions, $\mathbf{z}$. This objective directly addresses the primary goal of PC-EMS to reduce unnecessary ED visits while ensuring appropriate patient care. We formulate $D(\by, \bz)$ in Section \ref{sec_second_stage}. Constraint set \eqref{eq:constraint1} enforces that no more than $N_k$ units of type $k \in K$ are allocated to stations. Constraint set \eqref{eq:constraint2} maintains the vehicle ordering, ensuring that stations are populated with vehicles in a sequential manner. Constraint set \eqref{eq:constraint3} requires a unit $k \in K$ to be located at $i \in I$ if an ambulance of type $k$ is allocated from that station. Constraint set \eqref{eq:constraint33} requires that each patient at $j$ with believed condition $\theta$ be assigned at least one unit in the first stage. The inequality, $\ge$, allows multiple initial units to be dispatched to a single patient, allowing multiple response. Constraint \eqref{eq:constraint4} is a generic requirement that can be used to represent additional system constraints. Constraint sets \eqref{eq:constraint5} and \eqref{eq:constraint6} enforce the binary constraints.

\subsection{Second-Stage Model} \label{sec_second_stage}
The second-stage model determines how EMS resources are utilized once a patient’s true condition is revealed at the scene, given an allocation strategy $\bz$ and initial assignment decisions $\by$. This stage determines the specific actions performed and treatment pathways utilized by EMS resources. The model maximizes the number of patients diverted from the emergency department, while ensuring that all patients receive appropriate care while maintaining ambulance availability.

We define $\hat{\Theta}$ as the set of patient conditions that may be observed in the second stage. Let the parameter $p_{\theta \ttheta} \in [0,1]$ denote the probability that a patient initially classified as having condition $\theta \in \Theta$ actually has condition $\ttheta \in \hat{\Theta}$. Since all patients have a true condition, parameters $p_{\theta \ttheta}$ satisfy:
\begin{align}
    \sum_{\ttheta \in \hat{\Theta}} p_{\theta \ttheta} = 1, \forall \theta \in \Theta
\end{align}
If patient conditions are perfectly classified and $\Theta = \hat{\Theta}$, we have $p_{\theta \ttheta} = 1$ when $\theta = \ttheta$ and $p_{\theta \ttheta} = 0$ otherwise. Since initial classification is unlikely to be perfect, we expect $p_{\theta \ttheta} \in (0,1)$. Section \ref{uncertain_params} discusses how we estimate these probabilities.

To capture the decisions made after a patient’s true condition is revealed, we define the set $M_{ijk\ttheta}$ as the actions a unit from $i \in I$ of type $k \in K$ can perform when responding to a patient at $j \in J$ with an actual condition $\ttheta \in \hat{\Theta}$. This set includes both treatment pathways that satisfy a patient’s needs and assisting another responding unit without directly providing a treatment pathway.

We require each unit dispatched before the patient's true condition is learned, indicated by $\by$, to perform action $m \in M_{ijk\ttheta}$. The binary decision $x_{ijk \theta \ttheta } = 1$ indicates that a unit of type $k \in K$ from $i \in I$ performs action $m \in M_{ijk\ttheta}$ on patients at $j \in J$ with believed condition $\theta \in \Theta$ when their true condition is realized to be $\ttheta \in \hat{\Theta}$. 

To ensure every patient receives necessary care, at least one unit must perform an action from the subset $M'_{ijk\ttheta} \subset M_{ijk\ttheta}$, which represents treatment pathways that meet a patient’s needs. This requirement can be satisfied by a primary unit dispatched before the condition $\ttheta$ is known or by a secondary assignment unit. When secondary assignment is used, we let the binary decision variable $\hat{y}_{ijk \theta \ttheta} = 1$ indicate that a unit of type $k \in K$ from $i \in I$ is dispatched when a patient at $j \in J$ with believed condition $\theta \in \Theta$ is revealed to have actual condition $\hat{\theta} \in \hat{\Theta}$. We use the binary decision variable $\hat{x}_{ijk m \theta \ttheta} = 1$ to indicate the action $m \in M'_{ijk\ttheta}$ performed by this unit. 

Lastly, we define $M''_{ijk\ttheta} \subset M_{ijk\ttheta}$ as the subset of actions that result in a patient being diverted from the ED. The parameter $\lambda_{j \theta}$ represents the arrival rate of patients at $j \in J$ with an initially believed condition $\theta \in \Theta$. The term $\lambda_{j \theta} p_{\theta \ttheta}$ gives the arrival rate of patients at $j$ with an initially believed condition $\theta$, and actual condition $\ttheta$. The second-stage model is formulated as follows:
\begin{align}
     D(\by, \bz ) =  \maximize_{\mathbf{x},\mathbf{\hat{x}},\mathbf{\hat{y}}} \text{ } & \sum_{j \in J} \sum_{\theta \in \Theta} \sum_{\ttheta \in \hat{\Theta}} \lambda_{j \theta} p_{\theta \ttheta} \left[  \sum_{i \in I} \sum_{k \in K} \sum_{m \in M''_{ijk\ttheta}} (x_{ijkm \theta \ttheta} + \hat{x}_{ijkm \theta \ttheta}) \right] 
    \label{eq:s_objective} \\
    \text{subject to: } 
    & \sum_{i \in I} \sum_{k \in K} \sum_{m \in M'_{ijk\ttheta}} (x_{ijkm \theta \ttheta} + \hat{x}_{ijkm \theta \ttheta}) = 1,\quad  \forall j \in J, \theta \in \Theta, \ttheta \in \hat{\Theta}, \label{eq:s_c1} \\ 
    & \hat{y}_{ijk \theta \ttheta} \le z_{ik1},\quad \forall i \in I, j \in J, k\in K, \theta \in \Theta, \ttheta \in \hat{\Theta}, \label{eq:s_c2}\\
    & \sum_{m \in M_{ijk\ttheta}} x_{ijkm \theta \ttheta }  =  y_{ijk\theta},\quad \forall i \in I, j \in J, k \in K, \theta \in \Theta,\ttheta \in \hat{\Theta}, \label{eq:s_c3}\\
    &  \sum_{m \in M'_{ijk\ttheta}} \hat{x}_{ijkm \theta \ttheta} =  \hat{y}_{ijk \theta \ttheta},\quad \forall i \in I, j \in J, k \in K, \theta \in \Theta,\ttheta \in \hat{\Theta}, \label{eq:s_c4}\\
    & (\bx, \bxh, \by) \in P_{\alpha}(\by,\bz), \label{eq:s_c5}\\
     & x_{ijkm\theta\ttheta},\hat{x}_{ijkm\theta\ttheta} \in \{0,1\},  \forall i \in I, j \in J, k \in K, \theta \in \Theta,\ttheta \in \hat{\Theta}, m \in M_{ijk\hat{\Theta}}, \label{eq:s_c6}\\
     & \hat{y}_{ijk\theta \ttheta}\in \{0,1\}, \forall i \in I, j \in J, k \in K, \theta \in \Theta,\ttheta \in \hat{\Theta}. \label{eq:s_c7}
\end{align}
The objective function \eqref{eq:s_objective} maximizes the expected number of patients successfully diverted from the ED 
by adding up all actions from the set $M''_{ijk\tilde{\theta}}$. 
Constraint set \eqref{eq:s_c1} ensures that one unit performs an action from the subset $M'_{ijk\tilde{\theta}}$, which contains treatment pathways capable of meeting a patient's needs. Constraint \eqref{eq:s_c2} enforces that a secondary assignment decision $\hat{y}_{ijk \theta \tilde{\theta}}$ is only feasible if the first-stage allocation $\mathbf{z}$ allows for the deployment of unit type $k$ at station $i$. Constraint \eqref{eq:s_c3} ensures that any unit dispatched in the first stage must take exactly one action upon arrival, while constraint \eqref{eq:s_c4} guarantees that a secondary assigned unit, if activated, performs a valid treatment action. Constraint \eqref{eq:s_c5} stipulates that the assignment and action decisions must fall within a set $P_{\alpha}(\mathbf{y}, \mathbf{z})$, which represents the feasible region of dispatching decisions that maintain system reliability at level $\alpha$. This constraint ensures ambulance availability by limiting the traffic intensity at each station according to an $M/G/d/d$ queuing model, as formulated in Section \ref{sec_service_level}. Finally, constraints \eqref{eq:s_c6}--\eqref{eq:s_c7} enforce binary restrictions on the decision variables.

It is important to note that the feasibility of this two-stage model depends on the interplay between the constraints in the first stage (specifically those in $\mathcal{C}$), the second-stage constraints, and the availability requirements imposed by $P_{\alpha}(\mathbf{y}, \mathbf{z})$. Infeasibility can arise in several scenarios: if the parameter $\alpha$ is set too low (requiring very high ambulance availability), if the fleet size $N_k$ is insufficient, or if requirements in $\mathcal{C}$ are too stringent. When implementing this model, decision-makers should carefully balance competing constraints, potentially relaxing certain requirements when necessary to ensure feasibility while still meeting essential operational goals.

\subsection{Service-Level Guarantee} \label{sec_service_level}
Ensuring the availability of ambulances when needed is a key challenge in PC-EMS. Unlike traditional EMS models, where most patients follow a single transport pathway to the ED, PC-EMS introduces greater variability in resource utilization due to alternative treatment options and secondary assignments. To address this, we develop the feasibility set $P_{\alpha}(\by, \bz)$ in our second-stage model to enforce ambulance availability under system uncertainty. Specifically, this feasibility set ensures that when a new emergency call arrives, the probability of all ambulances at a station being unavailable does not exceed $\alpha$. This availability guarantee is constructed by incorporating constraints that limit the workload for each unit type at each station.

\subsubsection{Traffic Intensity and Blocking Probability.} 
To model ambulance availability, we represent each EMS station as an $M/G/d/d$ queue, where $M$ denotes Markovian arrivals (Poisson process), $G$ indicates general service time distribution, and the two $d$ parameters represent both the number of servers and maximum system capacity. This queueing model captures the random nature of emergency calls, accommodates variable service durations for different treatment pathways, and reflects that requests are lost when all ambulances are busy—all consistent with empirical EMS data \citep{kim_are_2014}.

Following approaches developed in previous work \citep{marianov_locationallocation_2002, boutilier_drone_2022}, we implement a chance constraint that limits the blocking probability at each station to at most $\alpha$. For each station $i$ and unit type $k$, we define the total traffic intensity as:
\begin{align} 
\rho_{ik} = \hat{\lambda}_{ik} \hat{s}_{ik}, 
\end{align} 
where $\hat{\lambda}_{ik}$ represents the average arrival rate of EMS requests at station $i$ for units of type $k$, and $\hat{s}_{ik}$ represents the corresponding average service time.

The blocking probability for an $M/G/d/d$ queue has a closed-form expression, allowing us to define a parameter $\rho_{\alpha d}$ as the maximum allowable traffic intensity that ensures the probability of all units being occupied does not exceed $\alpha$ \citep{boutilier_drone_2022}. The parameter $\rho_{\alpha d}$ increases monotonically with both $\alpha$ and $d$. Intuitively, this means that if we increase the likelihood of unit unavailability (higher \( \alpha \)) or increase the number of EMS units at a station (larger \( d \)), the system can accommodate a greater traffic load while maintaining the desired service level. Formally, if there are \( d \) units of type \( k \) at station \( i \), we ensure the probability of all units being occupied does not exceed $\alpha$ under $M/G/d/d$ assumptions by imposing the constraint:
\begin{align}
    \hat{\lambda}_{ik} \hat{s}_{ik} \leq \rho_{\alpha d} \label{eq:c_1}.
\end{align}

Next, we define the terms \( \hat{\lambda}_{ik} \) and \( \hat{s}_{ik} \) using the second-stage decision variables. Let \( \tau_{ijk \theta \ttheta} \) and \( \hat{\tau}_{ijk \theta \ttheta} \) represent the average service time a unit of type \( k \) from station \( i \) spends when assisting a patient at location \( j \) with a believed condition \( \theta \) and an actual condition \( \ttheta \), for initial and secondary assignments, respectively. We then define the average request rates (\( \hat{\lambda}_{ik} \)) and service times (\( \hat{s}_{ik} \)) based on actual demand and dispatch strategies as follows:
\begin{align}
    \hat{\lambda}_{ik} &= \sum_{j \in J} \sum_{\theta \in \Theta} \lambda_{j \theta} 
    ( y_{ijk \theta} + \sum_{\ttheta \in \hat{\Theta}} p_{\theta \ttheta} \hat{y}_{ijk \theta \ttheta} ), \\
    \hat{s}_{ik} &= \frac{\sum_{j \in J} \sum_{\theta \in \Theta} \sum_{\ttheta \in \hat{\Theta}} 
    \lambda_{j \theta} p_{\theta \ttheta} (\tau_{ijk \theta \ttheta} + \hat{\tau}_{ijk \theta \ttheta})}
    {\sum_{j \in J} \sum_{\theta \in \Theta} \lambda_{j \theta} 
    \left( y_{ijk \theta} + \sum_{\ttheta \in \hat{\Theta}} p_{\theta \ttheta} \hat{y}_{ijk \theta \ttheta} \right)}.
\end{align}
Given these definitions, constraint \eqref{eq:c_1} is equivalent to:
\begin{align}
   \sum_{j \in J} \sum_{\theta \in \Theta} \sum_{\ttheta \in \hat{\Theta}} 
   \lambda_{j \theta} p_{\theta \ttheta} (\tau_{ijk \theta \ttheta} + \hat{\tau}_{ijk \theta \ttheta}) 
   \leq \sum_{n = 1}^{N_k} (\rho_{\alpha n} - \rho_{\alpha (n-1)}) z_{ikn}, \quad \forall i \in I, k \in K. \label{eq:c_2}
\end{align}
The left-hand side represents the total expected service time required for EMS units at station \( i \) of type \( k \), based on realized patient conditions and treatment pathways. The right-hand side is equivalent to $\rho_{\alpha d}$ through the ordering constraint set \eqref{eq:constraint2} on \( z_{ikn} \),
where \( d \) denotes the number of units of type \( k \) at station \( i \), as indicated by \( \bz \).

\subsubsection{Service Time Computation.}
To quantify ambulance availability in our service-level constraints, we must accurately calculate the time units spend serving patients, including direct service, traveo times, and waiting periods. Since the average service times depend on the actions performed by EMS units, we define $\tau_{ijk\theta\ttheta}$ and $\hat{\tau}_{ijk\theta\ttheta}$ as decision variables that are a function of $\bx$, $\by$, $\bxh$ and $\byh$.  Let the parameter \( q_{ijkm\ttheta} \) represent the average service time required for a unit from station $i$ of type $k$ to perform action $m \in M_{ijk\ttheta}$ when treating a patient at location $j$ with an actual condition $\ttheta$. Additionally, let $r_{ij}$ represent the travel time from station $i$ to location $j$. 

The average time a unit of type $k$ from station $i$ spends assisting a patient at location $j$ with a believed condition $\theta$ and an actual condition 
$\ttheta$, denoted as $\tau_{ijk\theta\ttheta}$, consists of the average time required for actions performed during the response, plus the average delay resulting from waiting for a secondary assignment. The decision variable $\tau_{ijk\theta\ttheta}$ is defined as follows:
\begin{align}
    & \tau_{ijk\theta \ttheta} = \max \left\{ \sum_{m \in M_{ijk\ttheta}} q_{ijkm\ttheta} x_{ijkm \theta \ttheta} +  \sum_{i' \in I}  \sum_{k' \in K} \left[ r_{i'j} ( \hat{y}_{i'jk'\theta \ttheta} + y_{ijk\theta} - 1)\right],0 \right\}. \label{tau_1}
\end{align}

In instances where a unit is dispatched (i.e., $y_{ijk\theta} = 1$), $\sum_{m \in M_{ijk\ttheta}} x_{ijkm \theta \ttheta} = 1$ from \eqref{eq:s_c3}. In these cases, $\tau_{ijk\theta \ttheta}$ equals the sum of the average service time for the undertaken actions ($\sum_{m \in M_{ijk\ttheta}} q_{ijkm\ttheta} x_{ijkm \theta \ttheta}$) plus any time spent awaiting the arrival of a secondary assignment unit ($\sum_{i' \in I} \sum_{k' \in K} [ r_{i'j} (\hat{y}_{i'jk'\theta \ttheta})]$). Conversely, if no unit is dispatched (i.e., $y_{ijk\theta} = 0$), then $\sum_{m \in M_{ijk\ttheta}} x_{ijkm \theta \ttheta} = 0$ as mandated by \eqref{eq:s_c3}, resulting in $\tau_{ijk\theta\ttheta} = 0$.

The decision variable $\hat{\tau}_{ijk \theta \ttheta}$ represents the average service time of a unit from station $i$ of type $k$ that responds to a request at location $j$ that was believed to be condition $\theta$ and revealed to be condition $\ttheta$ as a secondary assignment unit. We formulate $\hat{\tau}_{ijk \theta \ttheta}$ as:
\begin{align}
     & \hat{\tau}_{ijk\theta \ttheta} = \sum_{m \in M_{ijk\ttheta}'} q_{ijkm\ttheta} \hat{x}_{ijkm \theta \ttheta}.  \label{tau_2}
\end{align}

We now formulate the complete service-level chance-constraint $P_\alpha(\by, \bz)$ as the set of feasible solutions $(\bx,\bxh,\byh)$ that satisfy the following linear constraints:
\begin{align}
   & P_\alpha(\by, \bz) = \Big\{ (\bx,\bxh,\byh) : \\
    &\sum_{j \in J} \sum_{\theta \in \Theta} \sum_{\ttheta \in \hat{\Theta}} 
 \lambda_{j \theta} p_{\theta \ttheta} (\tau_{ijk \theta \ttheta} + \hat{\tau}_{ijk \theta \ttheta}) \le  \sum_{n = 1}^{N_k} (\rho_{\alpha n} - \rho_{\alpha (n-1)}) z_{ikn}, \quad \forall i \in I , k \in K,  \label{p_queue1}\\
        &\sum_{m \in M_{ijk\ttheta}} q_{ijkm\ttheta} x_{ijkm \theta \ttheta} +  \sum_{i' \in I}  \sum_{k' \in K} \left[ r_{i'j} ( \hat{y}_{i'jk'\theta \ttheta} + y_{ijk\theta} - 1)\right] \le \tau_{ijk\theta \ttheta}, \forall i \in I, j \in J, k \in K, \theta \in \Theta, \ttheta \in \hat{\Theta}, \label{p_queue2}\\
    &\sum_{m \in M_{ijk\ttheta}'} q_{ijkm\ttheta} \hat{x}_{ijkm \theta \ttheta}  = \hat{\tau}_{ijk\theta \ttheta}, \quad \forall i \in I, j \in J, k \in K, \theta \in \Theta, \ttheta \in \hat{\Theta},  \label{p_queue3}\\
     &0\le \tau_{ijk\theta \ttheta},   \quad \forall i \in I, j \in J, k \in K, \theta \in \Theta, \ttheta \in \hat{\Theta}
       \Big\}. \label{p_queue4}
\end{align}

Constraint set \eqref{p_queue1} enforces the condition in \eqref{eq:c_2}, ensuring that the total expected service time does not exceed the allowable capacity for each station and unit type. Constraint sets \eqref{p_queue2} and \eqref{p_queue4} provide a linear formulation of equation \eqref{tau_1}, which accounts for the average service time incurred by both initial unit actions and delays from waiting for secondary assignments. Constraint set \eqref{p_queue3} corresponds directly to equation \eqref{tau_2}, capturing the average service time of any secondary-dispatched unit. Since these constraints are linear, they can be directly integrated into the second-stage integer programming model, ensuring the desired service-level guarantees.

\subsection{Model Discussion and Solution Approach}\label{sec:solution_approach}
The PC-EMS model presented above provides an accurate representation of a patient-centered emergency medical services system, but its size and complexity can make it computationally intensive to solve directly. To address this challenge, we implement a progressive solution approach that begins with a simplified version of the problem.

First, we limit each patient to a single response unit using the constraint:
\begin{align}
    \sum_{i \in I} \sum_{k \in K} y_{ijk\theta} \leq 1,
    \quad \forall\, j \in J,\;\theta \in \Theta.
    \label{eq:single_response}
\end{align}
This eliminates multiple response capabilities by ensuring that no more than one unit is initially dispatched to each patient. Second, we prohibit secondary assignments with:
\begin{align}
    \hat{y}_{ijk \theta \ttheta} = 0,
    \quad \forall\, i \in I,\, k \in K,\, j \in J,\, \theta \in \Theta,\, \ttheta \in \hat{\Theta}.
    \label{eq:no_secondary}
\end{align}

To solve the PC-EMS allocation model efficiently, we implement a progressive approach. We first solve a simplified version with constraints \eqref{eq:single_response}–\eqref{eq:no_secondary}, which limit each patient to a single response unit and prohibit secondary assignments. This represents standard dispatching policies in many EMS agencies while reducing computational complexity. The solution serves as a warm start when we subsequently relax these constraints to allow more complex dispatching strategies in the full model that combines the first, second stage, and service level-constraints (complete formulation in Online Appendix \ref{app_full_model}).

\section{Estimating the Uncertainty Parameter} \label{uncertain_params}
The probability parameter $p_{\theta \ttheta}$ plays a crucial role in linking the first and second stages of our PC-EMS allocation model. It determines the proportion of patients eligible for each alternative treatment pathway and the accuracy of dispatcher predictions before arrival. While prior studies suggest that phone screening can effectively distinguish patients requiring ED care \citep{lake_quality_2017}, dispatcher assessments in PC-EMS systems inherently contain uncertainty that must be quantified to develop effective allocation strategies.

To address this, we develop a machine learning approach to estimate $p_{\theta \ttheta}$. First, we construct a two-stage \textit{Diversion Eligibility Classification (DEC)} model to identify which patients are likely to be diverted from the ED and determine whether they would receive TIP or be transported to an AD. Then, we develop a \textit{Dispatcher Screening Prediction (DSP)} predictive model, only using information available during phone screening, that approximates how accurately EMS dispatchers can identify these patients before the onsite assessment. By comparing these predictions to patient eligibility, we quantify the uncertainty in dispatcher assessments.

This machine learning framework is not intended to be deployed for real-time clinical triage. Rather, it provides structured estimates of diversion potential and dispatcher accuracy to support system design and policy evaluation in PC-EMS.

\subsection{Data Source} \label{sec:data_source}
Data was obtained from the National EMS Information System (NEMSIS), a national database of EMS data from the United States \citep{national_highway_traffic_safety_administration_national_2021}. We utilize ground-based ambulances that treat and transport a single patient following an unscheduled emergent (911) request from September to December 2021. We perform separate analyses for urban and rural EMS systems. In this section, we summarize the urban analysis, while results from the rural analysis are provided in Online Appendix \ref{app_rural_analysis}. 

After filtering for completeness, the urban dataset consists of 515,806 unique EMS responses which we utilize for training the predictive model (257,903 responses), testing the predictive model (128,951 responses) and simulating our system (128,952 responses). Features include demographics, phone-screening complaints, procedures, and medications. To improve model stability, we removed highly correlated features ($|$correlation coeff$| > 0.85$) and encoded multi-label features as binary indicators, expanding the dataset from 24 to 407 features.

\subsection{Diversion Eligibility Classification (DEC)}
The first step in estimating $p_{\theta \ttheta}$ is to estimate which patients could be safely diverted from the ED in a PC-EMS system. While most EMS patients are transported to the ED in existing systems, some receive alternative treatments, either by declining transport or through protocol-permitted diversion. To estimate eligibility for ED diversion in a PC-EMS system, we utilize a machine-learning framework to learn the characteristics of patients that are currently diverted.

The first component of our DEC model predicts whether a patient can be diverted from the ED. Patients are labeled as diverted if they were treated using TIP or transported to an AD and not diverted if they were transported to an ED. To make this prediction, the model utilizes features that capture information about the patient and the onsite care the patient received. These independent variables include patient demographics (e.g. age, sex) to account for differences in care needs, clinical indicators (e.g. chief complaint, reported severity) to reflect medical urgency, EMS interventions (e.g. administered medications, performed procedures) to distinguish patients requiring higher levels of care, and dispatch characteristics (reason for EMS activation).

\begin{figure}[ht!]
\centering
\includegraphics[width=\textwidth]{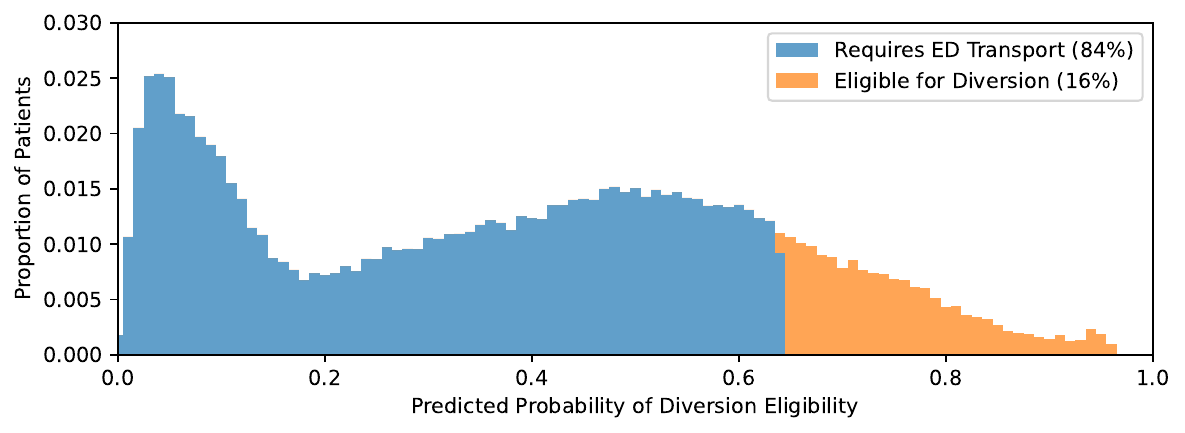}
\caption{Distribution of predicted probabilities used to estimate diversion eligibility. A classification threshold of 0.65 was selected to reflect prior estimates that 16\% of EMS patients may be diverted.}
\label{fig:curve1}
\end{figure}

Several classification algorithms were trained using three-fold cross-validation with balanced weighting to mitigate class imbalances on a subset of 100,000 EMS responses. A random forest classifier was selected for its ability to capture complex, non-linear relationships in the data. This model was applied to all 515,806 unique EMS responses. Since PC-EMS aims to expand diversion beyond current rates, the classification threshold was adjusted to include additional patients who share similar characteristics with those already diverted. Prior research suggests that approximately 16\% of EMS patients could be safely diverted \citep{department_of_health_and_human_services_innovation_2013}, so we set the classification threshold to 0.65 to reflect this expected proportion. The distribution of predicted probabilities is shown in Figure~\ref{fig:curve1}, the bimodal distribution suggests a separation between patients that may be suitable for diversion (right peak) and those requiring ED care (left peak).

The second component of the DEC classifies patients identified as diversion-eligible into either TIP or AD pathways. Patients who received on-site treatment without transport are labeled as TIP, while those transported to a non-ED facility are classified as AD. This model uses the same independent variables as the first component. Since the dataset of already diverted patients was much smaller than the full EMS dataset, training data was balanced to include an equal number of AD and TIP cases to prevent overfitting to the more common TIP class. As a result, the model was trained on a smaller dataset of 3,042 samples. A random forest classifier was selected for its ability to capture complex, non-linear relationships in the data. 

Once trained, this model is applied to the subset of patients identified as diversion-eligible by the first component of the DEC model. Since no prior information exists about the proportion of patients who would receive TIP versus AD in a PC-EMS system, we used a 0.5 classification threshold to assign patients based on learned patterns. Among those classified as diversion-eligible, the vast majority (94.6\%) were assigned to TIP, while 5.4\% were assigned to AD. While this distribution is imbalanced, it aligns with expectations. Patients requesting emergency services are more likely to require basic urgent care than specialized services available at an AD. Together, the DEC model provides labels indicating whether a patient is eligible for diversion and, if so, which alternative care pathway is most appropriate. These labels serve as a benchmark for evaluating how accurately EMS dispatchers can predict diversion eligibility, which we assess in the following section.

\subsection{Dispatcher Screening Prediction (DSP)}
After labeling patients likely to be diverted to TIP or an AD in a PC-EMS system, we develop the DSP model to assess how accurately EMS dispatchers can predict these labels using only information available during phone screening. Unlike the DEC, which incorporates on-site assessment details, the DSP model is restricted to pre-arrival data collected during the initial emergency call. By comparing its predictions to the labels generated by the DEC, we estimate how effectively human dispatchers could screen patients for diversion in a PC-EMS system. This model does not serve as an automated classification tool, instead approximating human decision-making under a structured dispatch process.

The DSP model predicts whether a patient is diversion-eligible based solely on pre-arrival information. Patients are labeled as diverted if the DEC model assigned them to TIP or AD and as non-diverted if they were transported to an ED. Independent variables include only phone-screening features available at the time of dispatch, such as patient demographics (e.g. age, sex), caller-reported symptoms, and the reason for EMS activation. These inputs reflect the limited information dispatchers use when making real-time diversion decisions.

\begin{figure}[ht!]
\centering
\includegraphics[width=\linewidth]{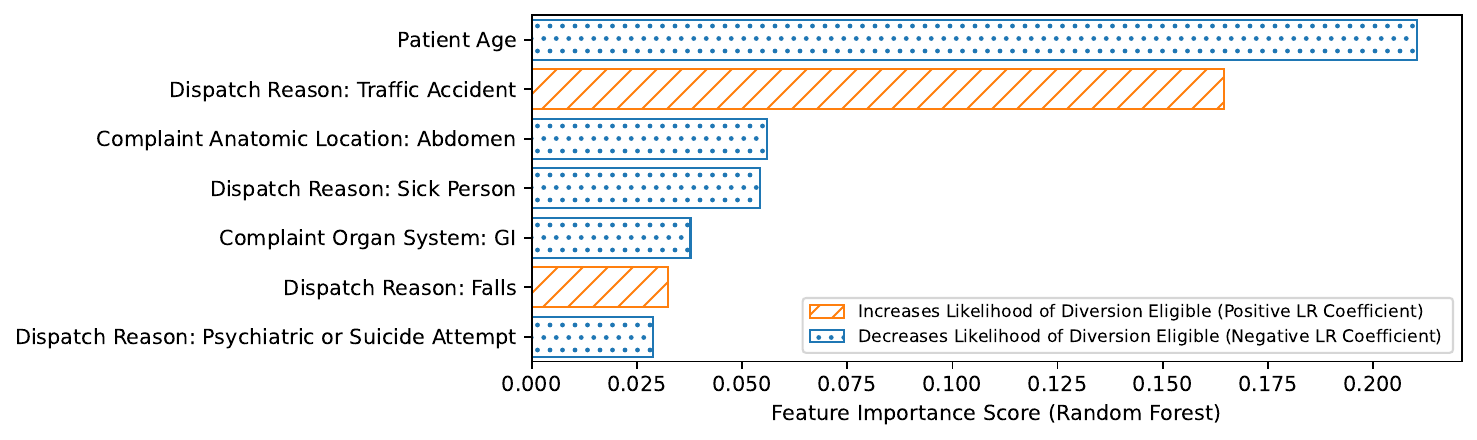}\caption{Bars indicate feature importance from a random forest model trained on phone-screening data. Hatching denotes the direction of association based on the corresponding coefficient in a logistic regression model.}
\label{fig:feature_importance}
\end{figure}

While a multinomial classifier could predict specific treatment pathways, we opted for a binary classification approach (diversion vs. no diversion) to reduce model complexity and improve computational efficiency in the optimization model. A random forest classifier was selected due to its alignment with EMS decision-making processes, where dispatchers rely on structured decision trees and screening protocols. The model was trained on 257,903 EMS responses collected between September and October 2021, using three-fold cross-validation with balanced weighting to mitigate class imbalances. Figure \ref{fig:feature_importance} presents feature importance scores from the random forest model. Notably, dispatches for specific events like falls and accidents are strongly associated with diversion eligibility, potentially because these incidents often require assessment but may not always necessitate ED care.

Once trained, the DSP model was applied to the test dataset to estimate dispatcher performance in predicting diversion-eligible patients. By comparing these predictions with the labels from the DEC, we quantify dispatcher uncertainty for the optimization framework.

Table \ref{tab:screening_accuracy} summarizes the probability parameters $p_{\theta \ttheta}$ used in our experiments. The sets $\Theta$ and $\hat{\Theta}$ correspond to treatment pathways, but they could be refined into more specific conditions. However, increasing the number of conditions at screening adds complexity for dispatchers, and expanding either set increases the computational cost of the optimization model. To balance feasibility and model accuracy, we maintain the classification at the treatment pathway level.

\begin{table}[ht!]
\centering
\renewcommand{\arraystretch}{.7} 
\caption{Estimated probability of actual treatment needs ($\hat{\Theta}$) conditioned on screening predictions ($\Theta$). These probabilities define the parameter $p_{\theta \ttheta}$ in our optimization model.}
\label{tab:screening_accuracy}
\begin{tabular}{|C{6cm}|C{3.1cm}||C{1.6cm}|C{1.6cm}|C{1.6cm}|}
\hline
\multirow{2}{*}{\shortstack[c]{\textbf{Believed Condition ($\Theta$)}}} & 
\multirow{2}{*}{\raisebox{-.25\height}{\shortstack[c]{\textbf{\% of Patients}\\ \textbf{Screened}}}} & \multicolumn{3}{c|}{\textbf{Actual Condition ($\hat{\Theta}$)}} \\ 
\cline{3-5}
& & \textbf{ED} & \textbf{AD} & \textbf{TIP} \\
\hline
Likely Not Diversion Eligible (ED)  & 70.3\%  & 93.2\%  & 0.4\%  & 6.3\%  \\
\hline
Likely Diversion Eligible (AD/TIP)  & 29.7\%  & 62.7\%  & 1.9\%  & 35.4\%  \\
\hline
\end{tabular}
\end{table}

\section{Case Study Set-up}\label{sec:case_study}
To evaluate our PC-EMS allocation model, we designed experiments across four distinct EMS environments with varying population densities, geographic scales, and demand distributions. 

\subsection{Data Sources and Environment Construction}
We built our test environment using three datasets. Geographic demand patterns are based on simulated cardiac arrest data developed by \cite{boutilier_drone_2022}, which capture realistic spatial variations across regions in Ontario, Canada. Patient-level data from the NEMSIS database (see Section~\ref{sec:data_source}) inform estimates of service times for different treatment pathways. Facility location data from Canadian health services records are used to identify hospitals with EDs and clinics serving as ADs \citep{statistics_canada_open_2020}.

\subsection{Geographic Regions}\label{sec:geo_regions}
We analyze four distinct EMS environments to identify how optimal PC-EMS strategies vary with regional characteristics. The urban system (Hamilton, ON) features a compact operating environment with high call density and short travel times, where resource availability typically constrains service delivery more than geographic coverage. The rural system (Muskoka, ON) is characterized by sparse demand distributed across extensive distances, where travel times significantly impact service delivery. The mixed system (Halton, ON) combines concentrated urban call volumes with dispersed demand in peripheral areas, creating varying response dynamics within a single region. The suburban system (Simcoe County, ON) spans a broader area with multiple medium-density population clusters, presenting intermediate challenges between urban and rural contexts.

\begin{figure}[ht!]
    \centering
    \subfloat[Urban]{%
        \parbox[c]{0.45\textwidth}{\centering
        \includegraphics[height = 2.5in]{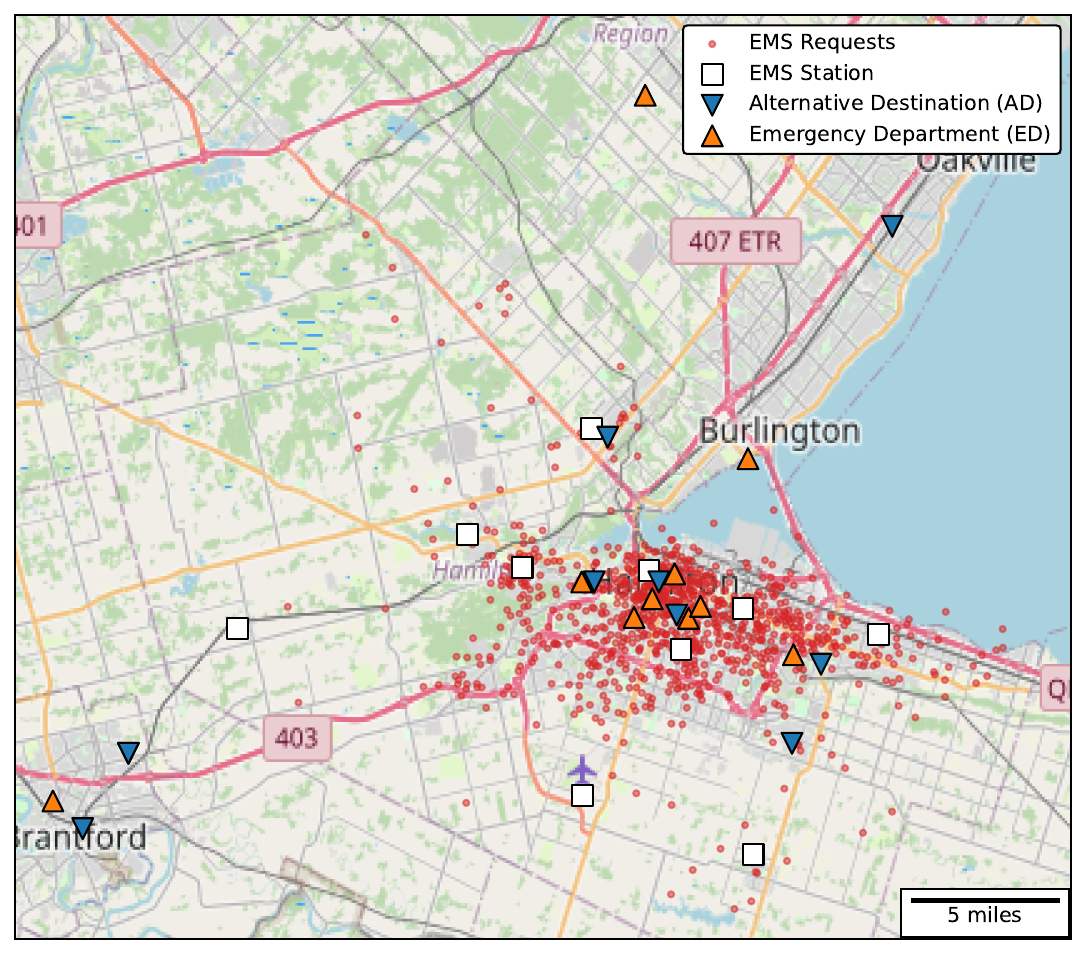}}%
        \label{fig:map_urban}%
    }%
    \hfill
    \subfloat[Rural]{%
        \parbox[c]{0.45\textwidth}{\centering
        \includegraphics[height = 2.5in]{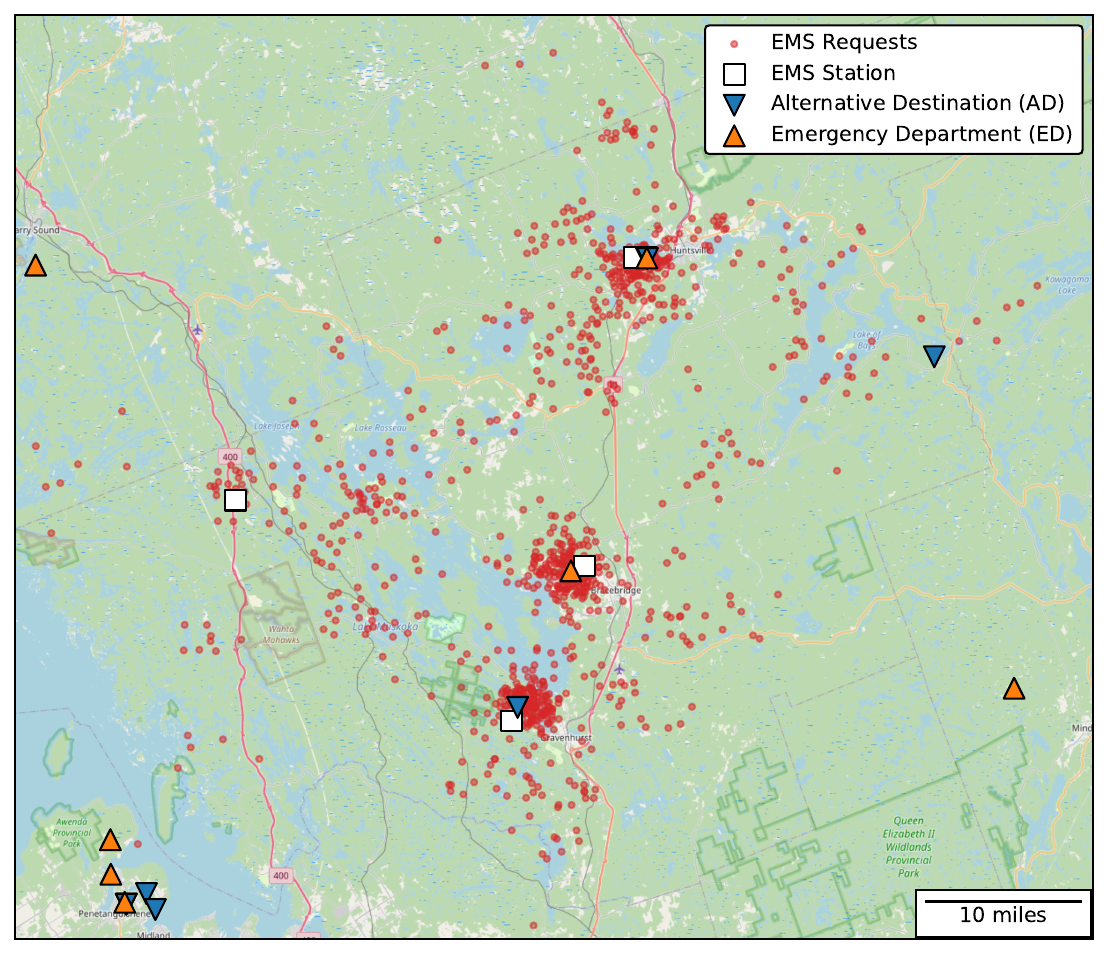}}%
        \label{fig:map_rural}%
    }
    
    \vspace{.2em} 
    
    \subfloat[Mixed]{%
        \parbox[c]{0.45\textwidth}{\centering
        \includegraphics[height = 2.7in]{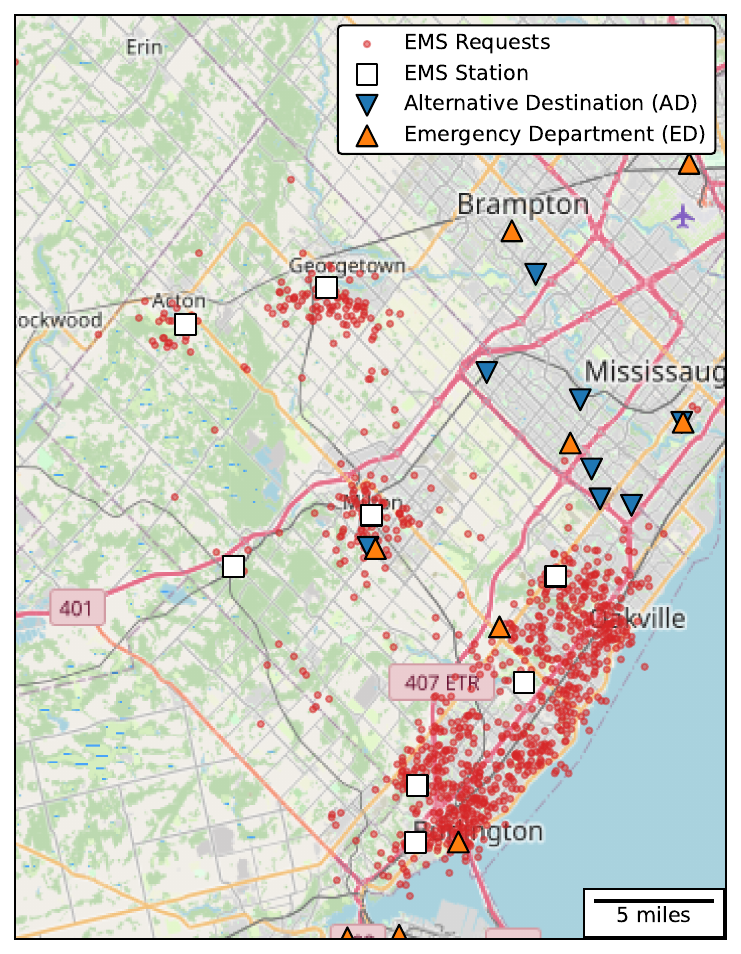}}%
        \label{fig:map_mixed}%
    }%
    \hfill
    \subfloat[Suburban]{%
        \parbox[c]{0.45\textwidth}{\centering
        \includegraphics[height = 2.7in]{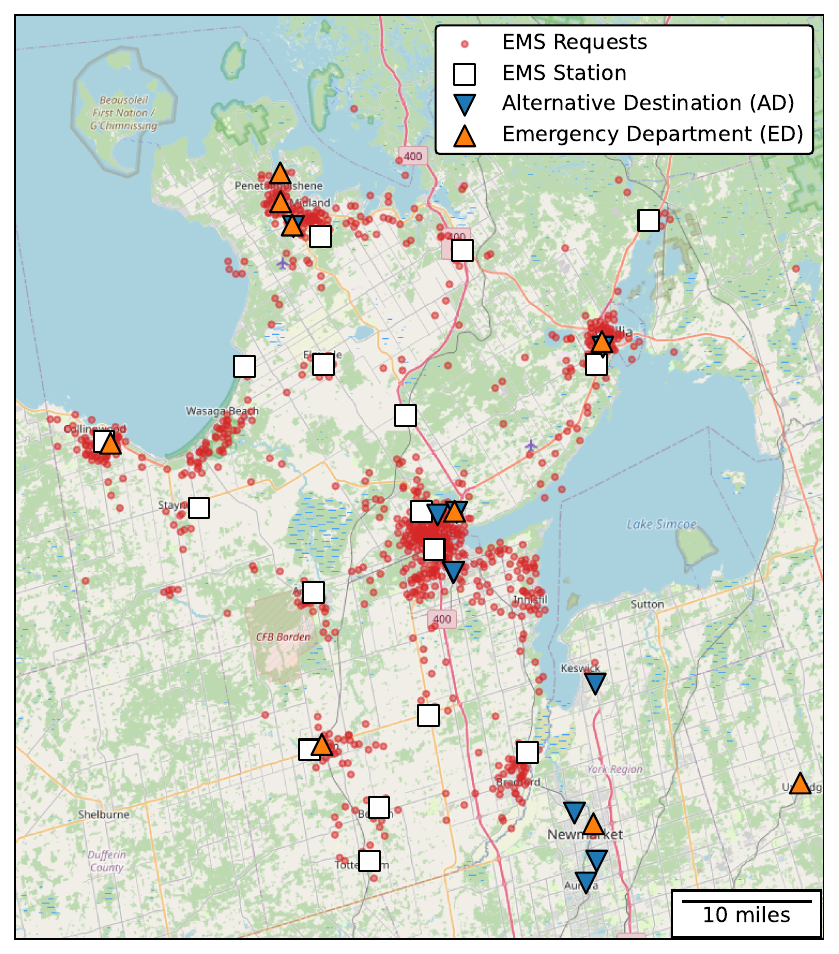}}%
        \label{fig:map_suburban}%
    }
    \caption{Geographic distribution of demand, ambulance stations, and healthcare facilities locations.}
    \label{fig:all_maps}
\end{figure}

Figure \ref{fig:all_maps} shows maps of the four geographic regions studied in our analysis. Each map displays the distribution of emergency call demand, ambulance station locations, hospitals with emergency departments, and urgent care and primary care offices that may serve as alternative destinations.

\subsection{System Parameters and Demand Modeling}\label{sec:system_params}
We constructed modeling parameters for each region as follows. Each geographic region was divided into a grid of equal-sized cells to form the set of demand nodes $J$. The arrival rate parameter $\lambda_{j\theta}$ was derived from cardiac arrest frequencies, scaled by the observed ratio of total EMS responses to cardiac arrests (157:1) from the NEMSIS dataset. We considered two conditions in the set $\Theta = \{\text{``likely diversion eligible"}, \text{``likely not diversion eligible"}\}$ and assign call probabilities according to the classification probabilities from Table~\ref{tab:screening_accuracy}.

The set of actual conditions $\hat{\Theta} = \{\text{``ED"}, \text{``AD"}, \text{``TIP"}\}$ represents the appropriate treatment option determined after on-scene assessment. The probability parameters $p_{\theta\hat{\theta}}$ capture the likelihood that a patient initially classified with condition $\theta$ actually has condition $\hat{\theta}$, based on our analysis in Section~\ref{uncertain_params}. The set of potential station locations $I$ was selected from existing EMS stations, identifying subsets that maintained maximum coverage while preserving at least 50\% of stations in each region. This approach balanced realistic deployment scenarios with computational requirements.

Table~\ref{tab:region_summary} summarizes the key characteristics of each region, including annual call volumes, facility counts, and geographic parameters. Travel times ($r_{ij}$) between locations were calculated using the travel time function from \citet{kolesar_determining_1975}, which estimates travel time based on Euclidean distance. Service time parameters ($q_{ijkm\hat{\theta}}$) were derived from NEMSIS data: ED transports required 49 minutes for on-scene care and transfer (excluding travel), AD transports averaged 43 minutes, and TIP interventions took approximately 45 minutes based on community paramedicine literature \citep{shah_improving_2018, dainty_home_2018}. Supporting units (i.e., those assisting primary responders without providing transport) were assigned 5 minutes of service time, a conservative estimate accounting for the lack of specific data on these interactions.

\begin{table}[ht!]
\centering
\renewcommand{\arraystretch}{.7}
\caption{Summary of PC-EMS Test Regions}
\begin{tabular}{|l|C{2.6cm}|C{2.6cm}|C{2.6cm}|C{2.6cm}|}
\hline
\textbf{Attribute} & \textbf{Urban} & \textbf{Mixed} & \textbf{Suburban} & \textbf{Rural} \\
\hline
\raisebox{-.5\height}{\shortstack[l]{Geographic data\\based on}}& \raisebox{-.5\height}{\shortstack[c]{City of \\Hamilton, ON} }& \raisebox{-.5\height}{\shortstack[c]{Region of \\Halton, ON} }& \raisebox{-.5\height}{\shortstack[c]{County of \\Simcoe, ON} }& \raisebox{-.5\height}{\shortstack[c]{District of \\Muskoka, ON}}  \\[7.5pt]
\hline
Annual EMS Calls & 97,185 & 55,826 & 69,139 & 9,523 \\
\hline
\# of EMS Stations & 10 & 8 & 17 & 4 \\
\hline
Total EMS Units & 25 & 24 & 38 & 8 \\
\hline
Demand Node Size (mi\textsuperscript{2}) & 1.5 & 1.5 & 2.0 & 5.0 \\
\hline
\# of Demand Nodes & 217 & 239 & 501 & 83 \\
\hline
\end{tabular}
\label{tab:region_summary}
\end{table}

To determine fleet sizes, we developed a coverage optimization model that identified the minimum number of total units needed to achieve both maximum coverage and 95\% ambulance availability ($\alpha = 0.05$) in each region. This approach is necessary because actual fleet sizes for the regions were unknown.

\subsection{PC-EMS Implementation}\label{sec:pcems_impl}
Our implementation includes two types of EMS units in set $K$: traditional units limited to standard ED transport protocols, and diversion-capable units equipped with advanced training or telemedicine capabilities that enable assessment for alternative treatment pathways. 
For the three conditions in $\hat{\Theta} = \{\text{``ED"}, \text{``AD"}, \text{``TIP"}\}$, we enforce a hierarchical treatment pathway structure where patients eligible for TIP can also be treated at an AD or ED, and patients eligible for AD can also be treated at an ED. This hierarchy defines our treatment pathways sets, with $M'_{ijk\hat{\theta}}$ representing valid treatment pathways and $M''_{ijk\hat{\theta}}$ representing treatment pathways resulting in ED diversion. Traditional units can only perform ED transport actions, while diversion-capable units can implement the full range of treatment options. The set $M_{ijk\hat{\theta}}$ also includes supporting actions for all unit types when assisting through multiple response or secondary assignment.

The first-stage model enforces minimum coverage constraints, ensuring that the implementation of PC-EMS maintains the same geographic accessibility standards as traditional EMS systems.

\subsection{Optimization and Simulation Framework}\label{sec:simulation}
To evaluate system performance under realistic conditions, our methodology employs a two-phase approach. In the optimization phase, we solve the PC-EMS allocation model for each region and configuration using the solution approach described in Section~\ref{sec:solution_approach}. These optimization problems were solved using Gurobi 9.5 on a high-performance computing cluster with 16 CPU cores and 256GB of RAM, with a maximum 24-hour solve time. Of the 714 problem instances solved across all geographic regions, 699 (97.9\%) achieved an absolute optimality gap below 0.01. The remaining 15 instances had optimality gaps between 0.01 and 0.24, all occurring in urban and suburban environments. We attribute these larger gaps to the higher density of potential station locations in these regions, creating more interchangeable solutions, slowing convergence.

In the simulation phase, we use the resulting allocations to parameterize a discrete-event simulation that captures the dynamic nature of emergency service operations. Emergency requests arrive according to non-homogeneous Poisson processes with rates varying by time of day and day of week, based on temporal patterns observed in the NEMSIS data. Call locations are sampled from the geographic demand distribution. Travel times and service durations follow exponential distributions, with means parameterized by the values estimated from our NEMSIS analysis.

The simulation implements the optimized dispatching strategy for each request. If recommended units are unavailable, the system follows a fallback protocol by dispatching the closest available ambulance and defaulting to ED transport. For each region and configuration, we performed 100 replications of 7-day simulations using a computing cluster with 8 CPU cores and 128GB of RAM. 

\section{Policy Experiments}
Our policy experiments aim to provide actionable insights for future PC-EMS implementations. We focus particularly on maximizing diversion rates, as they provide the strongest justification for PC-EMS from a health system perspective. 

\subsection{Fleet Composition for Effective Diversion} \label{sec:diversion}
Our model distinguishes between two types of EMS units. Traditional units restricted to ED transport and diversion-capable units that can implement alternative treatment pathways. Since upgrading an entire fleet would be costly, agencies must determine what proportion of units should have these capabilities. To identify the optimal fleet composition, we systematically varied the proportion of diversion-capable units while maintaining total fleet size, measuring the resulting diversion rates through simulation.

\begin{figure}[ht!]
    \centering
    \includegraphics[width=0.8\linewidth]{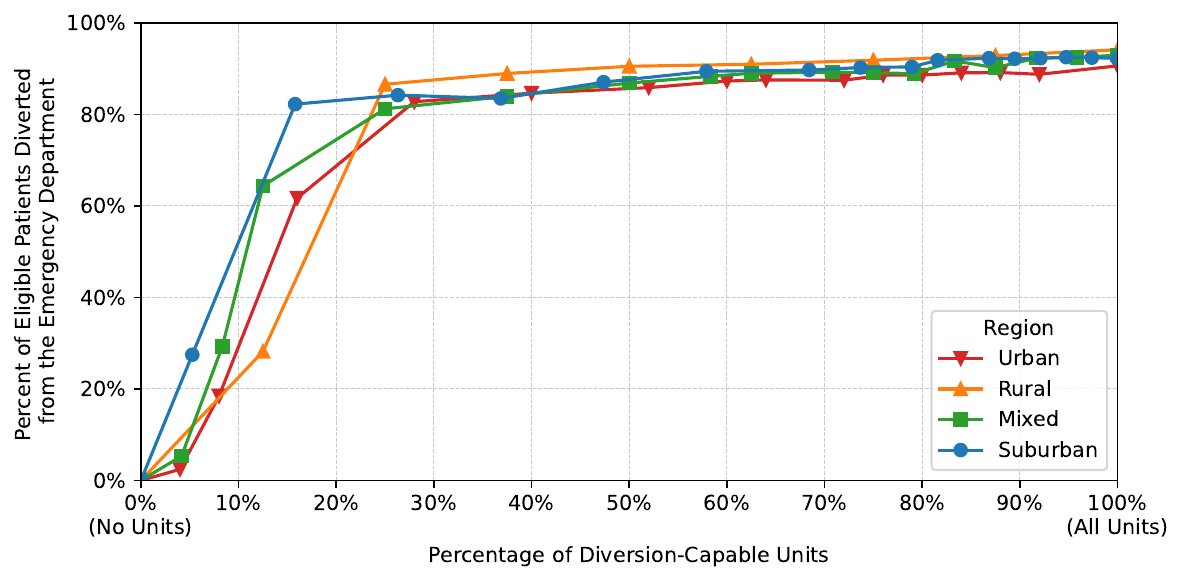}
   \caption{Percentage of potential diversions achieved versus percentage of diversion-capable units. Each region has a different total fleet size, but percentages allow for standardized comparison across settings.}
    \label{fig:diverion_fleet}
\end{figure}

Our analysis demonstrates that EMS agencies can achieve approximately 80\% of potential patient diversions by equipping only a fraction of their fleet with diversion capabilities. Specifically, in the urban environment, equipping 28\% of units (7 out of 25 units) with diversion capabilities achieved the 80\% threshold of potential diversions. In the rural environment, 25\% diversion-capable units (2 out of 8 units) was sufficient to reach this threshold. The mixed environment similarly required 25\% diversion-capable units (6 out of 24 units), while the suburban environment reached the 80\% threshold with only 16\% diversion-capable units (6 out of 38 units). Figure \ref{fig:diverion_fleet} illustrates this relationship between the percentage of diversion-capable units and the percentage of potential diversions achieved across all four geographic environments.

This finding has significant implications for the cost-effective implementation of PC-EMS. Rather than equipping the entire fleet with diversion capabilities, agencies can strategically equip a smaller fraction of their units using our optimization framework while still diverting the majority of eligible patients. This targeted approach reduces implementation costs while realizing the potential operational benefits of patient-centered care.

The sensitivity analysis performed assumes multiple response and secondary assignment are implemented. However, EMS agencies vary in their willingness to adopt such strategies. We explore the impact of these dispatching policies in the next section.

\subsection{Impact of Dispatching Strategies}
To evaluate the importance of multiple response and secondary assignment strategies in diverting patients, we compared three dispatching approaches: (1) \textit{Full Dispatch} with both multiple response and secondary assignment capabilities, (2) \textit{Multiple Response} without secondary assignment, and (3) \textit{Single Response} with one unit per patient only. \textit{Full Dispatch} is the full model and results presented in \ref{sec:diversion}.  For the \textit{Multiple Response} approach, we imposed constraint set \eqref{eq:no_secondary}, which eliminates the possibility of dispatching additional units after the patient's condition is revealed. For the \textit{Single Response} approach, we imposed constraint sets \eqref{eq:no_secondary} and \eqref{eq:single_response}, which restricts the system to assign at most one unit per patient, eliminating both secondary assignment and multiple response possibilities.

\begin{figure}[ht!]
    \centering
    \includegraphics[width=\linewidth]{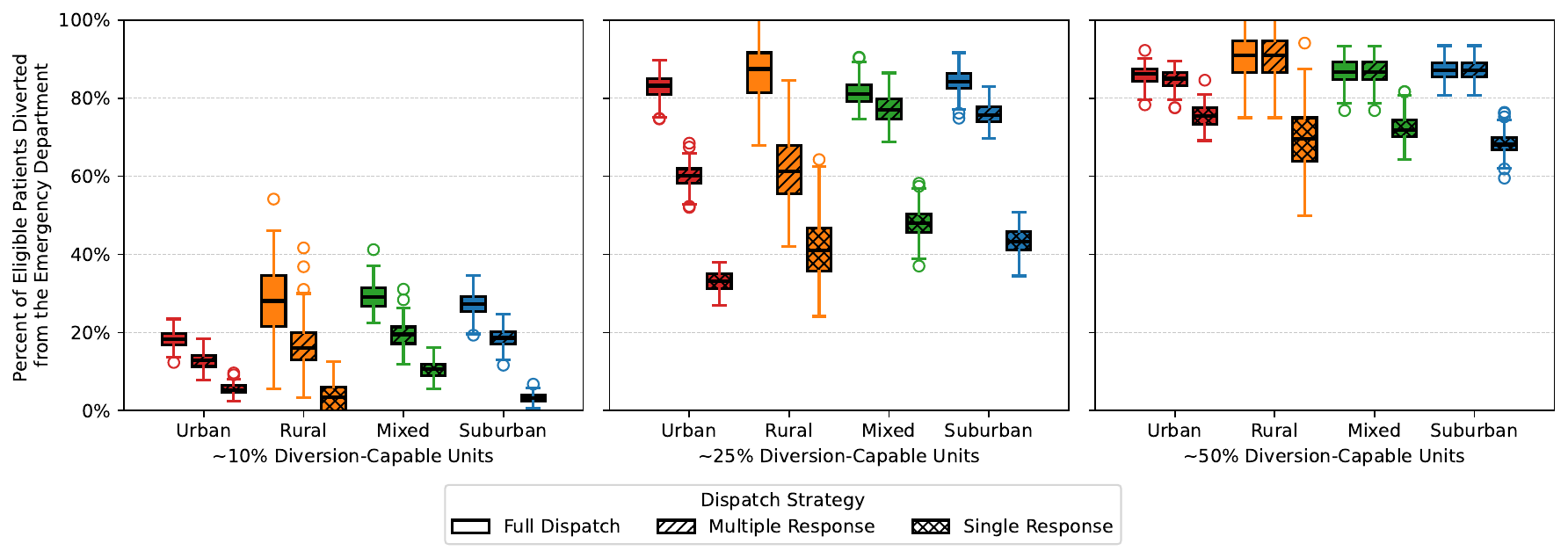}
    \caption{Impact of dispatching strategies on diversion rates with $\sim$10\%, 25\%, and 50\% diversion-capable units. The boxplots display the distribution of diversion rates observed across simulation replications.}
    \label{fig:dispatch_figure}
\end{figure}

Our analysis reveals that the \textit{Full Dispatch} strategy significantly enhances diversion rates when diversion-capable units comprise a smaller portion of the fleet. At approximately 10\% diversion-capable units, the \textit{Full Dispatch} model achieves approximately 3.4-8.6 times more diversions than the \textit{Single Response} approach, with the most dramatic improvement in suburban settings (27.5\% vs 3.2\%). At approximately 25\%, it achieves approximately 1.9-2.5 times more diversions, with the most substantial difference in urban areas (82.8\% vs 33.0\%). Figure \ref{fig:dispatch_figure} presents these impacts of dispatching strategies on diversion rates across all four environments for approximately 10\%, 25\%, and 50\% of diversion-equipped units.

The improvement of the \textit{Full Dispatch} strategy narrows as the proportion of diversion-capable units increases to approximately 50\%, with the improvement reducing to approximately 1.1-1.3 times across environments. \textit{Multiple Response} typically falls somewhere between the \textit{Full Dispatch} and \textit{Single Response} approaches, achieving approximately 60.1\% of potential diversions with 25\% diversion-capable units in urban areas. Notably, multiple response is typically easier for agencies to implement than secondary assignment, requiring less coordination and communication.

These findings suggest that EMS agencies should investigate procedures to support multiple response and secondary assignment dispatching strategies, particularly when budget constraints limit the number of diversion-capable units. Results for other all fleet compositions are provided in Online Appendix \ref{app_dispatch_sens}.

\subsection{Value of Patient Screening Accuracy} \label{sec:value_screening}
The PC-EMS allocation model explicitly accounts for uncertainty in patient classification through the parameter $p_{\theta \ttheta}$, which links first-stage dispatching decisions to second-stage treatment decisions. The ET3 program initially included funding for advanced phone screening capabilities (which was later cut), but the value of improved screening has not been quantified in the PC-EMS context \citep{centers_for_medicare__medicaid_services_emergency_2022}.

\begin{figure}[ht!]
    \centering
    \includegraphics[width=\linewidth]{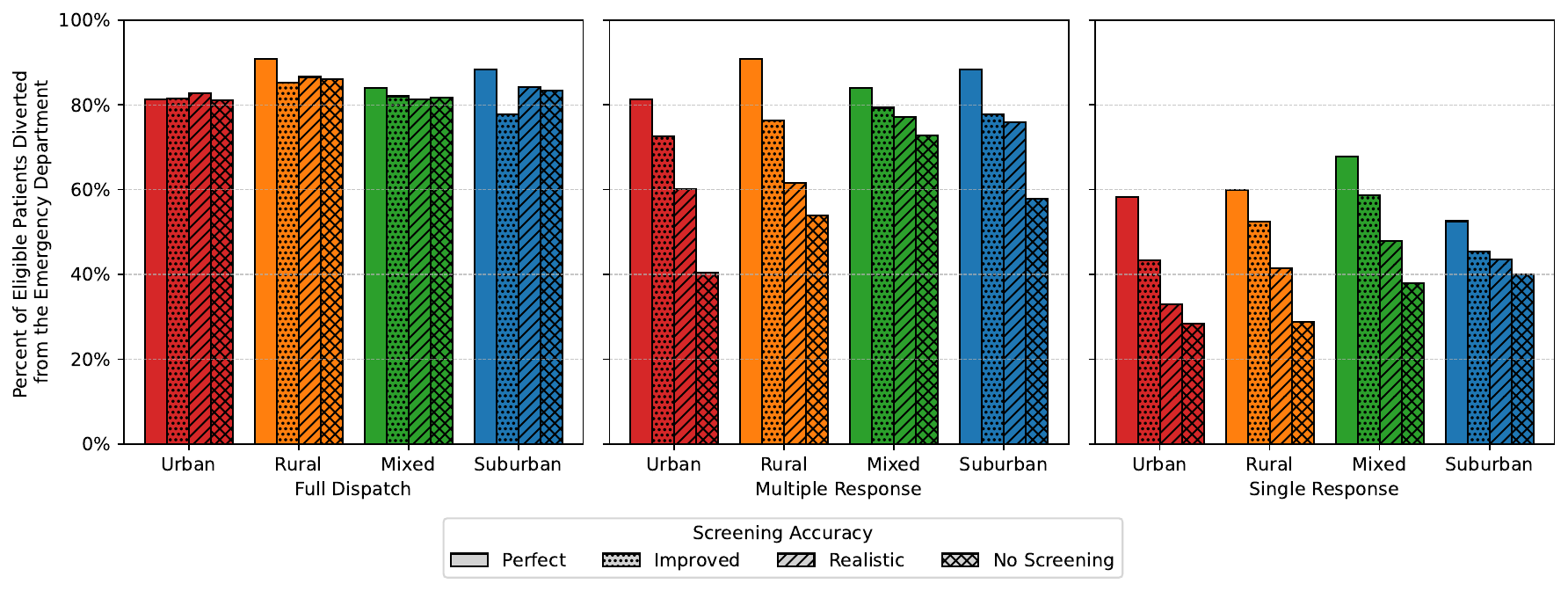}
    \caption{Effect of patient screening accuracy on diversion rates with $\sim 25\% $ diversion-capable fleet composition across different geographic environments and dispatching strategies.}
    \label{fig:screening_fig}
\end{figure}

To explore this relationship, we modified the accuracy parameters from our baseline model in several ways, examining the impact of: (1) perfect information, where patient conditions are known with certainty before dispatch; (2) no screening, where dispatching occurs without patient classification; (3) improved screening, which enhances the accuracy of the baseline parameters $p_{\theta \ttheta}$ by reducing misclassification rates by 50\% (i.e., increasing true positives and true negatives while decreasing false positives and false negatives identified in Section \ref{uncertain_params}); and (4) realistic screening, which uses the baseline parameters $p_{\theta \ttheta}$ we identified through our machine learning approach. Figure \ref{fig:screening_fig} demonstrates that screening accuracy significantly influences system performance when there is approximately 25\% diversion-capable units. Results for other fleet compositions are provided in Online Appendix \ref{app_screening_sens}.

We observe that screening accuracy has minimal effect when recourse through secondary assignment is allowed. With Full Dispatch, diversion rates remain consistently high across all screening scenarios (86.2\% with perfect information vs. 82.8\% with realistic screening in the urban environment). This makes intuitive sense, as secondary assignment provides opportunities to correct initial dispatching mistakes. Screening accuracy becomes more important in single response settings when the system cannot adjust the response, with diversion rates increasing from 33.0\% with realistic screening to 58.3\% with perfect information in urban environments.

These findings suggest that if agencies cannot implement multiple response and secondary assignment dispatching strategies, they should prioritize improving patient screening capabilities. Enhanced screening protocols and dispatcher training can significantly boost diversion rates in simpler dispatching systems, offering a potential alternative to implementing recourse mechanisms.

\subsection{Managerial Insights}
Our findings reveal that PC-EMS implementation faces a trade-off between fleet composition and coordination. While agencies can achieve approximately 80\% of potential diversions by equipping only 15-25\% of their fleet with advanced capabilities, this efficiency depends on dispatching strategies that may be challenging to implement in practice. Multiple response and secondary assignment protocols require significant coordination overhead, increased workload for providers, and modifications to existing workflows in many US EMS systems. Additionally, enhanced dispatcher screening capabilities represent another layer of implementation beyond fleet modifications.

When planning PC-EMS implementation, agencies should carefully assess their organizational readiness for operational change. Organizations that are hesitant to utilize multiple response and secondary assignment or are unable to invest in advanced patient screening may opt to equip a larger proportion of their fleet with diversion capabilities while maintaining simpler dispatching protocols—accepting higher equipment costs to minimize operational disruption. The optimal strategy depends on an agency's existing infrastructure, workforce capabilities, and tolerance for process change. Successful implementation ultimately requires a realistic assessment of both the potential benefits and the practical constraints involved in transforming emergency response systems.
\section{Conclusion}
This research addresses a critical challenge in emergency medical services by developing a comprehensive framework for PC-EMS allocation. The findings offer several critical insights for EMS agencies and policymakers. First, our analysis reveals that agencies can achieve up to 80\% of potential patient diversions by equipping only 20-25\% of their fleet with advanced capabilities. Second, multiple response and secondary assignment can significantly improve diversion rates, especially when only a limited number of units are diversion-capable. Finally, if agencies cannot implement these dispatching strategies, improved patient screening accuracy provides an alternative pathway to enhance diversion rates.

While the ET3 model demonstrated the potential of patient-centered emergency medical services, it was ultimately discontinued due to implementation challenges. Our research provides a structured approach to overcoming these obstacles, offering a data-driven methodology for transforming emergency medical response.

The principles developed in this study extend beyond emergency medical services. The proposed allocation model and uncertainty estimation framework can be adapted to other resource allocation problems, such as utility vehicle deployment during natural disasters or critical medical equipment distribution. Future research could explore more granular patient condition classifications, develop real-time decision support tools for EMS dispatchers, and investigate the long-term clinical and economic impacts of PC-EMS implementation. By balancing operational constraints with patient-centered care principles, we can create more flexible, efficient, and responsive emergency medical systems that better serve community health needs.

\ACKNOWLEDGMENT{
The authors acknowledge the National Emergency Medical Services Information System (NEMSIS) for providing the emergency medical services data used in this research. The content reproduced from the NEMSIS Database remains the property of the National Highway Traffic Safety Administration. The National Highway Traffic Safety Administration is not responsible for any claims arising from works based on the original Data, Text, Tables, or Figures.
\\
The authors acknowledge Statistics Canada for providing the Open Database of Healthcare Facilities (ODHF) used to identify healthcare facility locations in this research. The ODHF is released as open data under the Open Government License - Canada.
}

\bibliographystyle{custom.bst} 
\bibliography{main} 

\setcounter{page}{1}
\RUNAUTHOR{Stratman, Boutilier, and Albert}
\RUNTITLE{Online Appendix to "Ambulance Allocation for Patient-Centered Care"}
\section*{Online Appendix to ``Ambulance Allocation for Patient-Centered Care''}
\text{}
 \begin{APPENDICES}
 \numberwithin{equation}{section}
 \renewcommand{\theequation}{\thesection.\arabic{equation}}
        \numberwithin{figure}{section}
        \numberwithin{table}{section}
        \renewcommand{\thefigure}{\thesection.\arabic{figure}}
        \renewcommand{\thetable}{\thesection.\arabic{table}}
     \section{Complete Optimization Model Formulation}\label{app_full_model}
         The following mixed-integer linear programming model represents the full mathematical representation of the optimization problem. See the main text for interpretation and discussion of the model's components.
        \begin{align}
        & \maximize_{\mathbf{y},\mathbf{z},\mathbf{x},\mathbf{\hat{x}},\mathbf{\hat{y}},\mathbf{\tau},\mathbf{\hat{\tau}}} \text{ }  \sum_{j \in J} \sum_{\theta \in \Theta} \sum_{\ttheta \in \hat{\Theta}} \lambda_{j \theta} p_{\theta \ttheta} \left[ \sum_{i \in I} \sum_{k \in K} \sum_{m \in M''_{ijk\ttheta}} (x_{ijkm \theta \ttheta} + \hat{x}_{ijkm \theta \ttheta}) \right] \\
        &\text{subject to: } \\
        & \sum_{i \in I}\sum_{n = 1}^{N_k } z_{ikn} \le N_k, \quad \forall k \in K, \\
        & z_{ikn} \le z_{ik(n-1)}, \quad \forall i \in I, k \in K, n \in {2,...,N_k}, \\
        & y_{ijk\theta } \le z_{ik1}, \quad \forall i \in I, j \in J, k \in K, \theta \in \Theta, \\
            & \sum_{i \in I} \sum_{k \in K} y_{ijk\theta } \ge 1, \quad \forall j \in J, \theta \in \Theta,  \\ 
        & \sum_{i \in I} \sum_{k \in K} \sum_{m \in M'_{ijk\ttheta}} (x_{ijkm \theta \ttheta} + \hat{x}_{ijkm \theta \ttheta}) = 1, \quad \forall j \in J, \theta \in \Theta, \ttheta \in \hat{\Theta},  \\
        & \hat{y}_{ijk \theta \ttheta} \le z_{ik1}, \quad \forall i \in I, j \in J, k\in K, \theta \in \Theta, \ttheta \in \hat{\Theta}, \\
        & \sum_{m \in M_{ijk\ttheta}} x_{ijkm \theta \ttheta } = y_{ijk\theta}, \quad \forall i \in I, j \in J, k \in K, \theta \in \Theta,\ttheta \in \hat{\Theta}, \\
        & \sum_{m \in M'_{ijk\ttheta}} \hat{x}_{ijkm \theta \ttheta} = \hat{y}_{ijk \theta \ttheta}, \quad \forall i \in I, j \in J, k \in K, \theta \in \Theta,\ttheta \in \hat{\Theta}, \\
        & \sum_{j \in J} \sum_{\theta \in \Theta} \sum_{\ttheta \in \hat{\Theta}} \lambda_{j \theta} p_{\theta \ttheta} (\tau_{ijk \theta \ttheta} + \hat{\tau}_{ijk \theta \ttheta}) \le \sum_{n = 1}^{N_k} (\rho_{\alpha n} - \rho_{\alpha (n-1)}) z_{ikn}, \quad \forall i \in I , k \in K, \label{queue1}\\
        & \sum_{m \in M_{ijk\ttheta}} q_{ijkm\ttheta} x_{ijkm \theta \ttheta} + \sum_{i' \in I} \sum_{k' \in K} \left[ r_{i'j} ( \hat{y}_{i'jk'\theta \ttheta} + y_{ijk\theta} - 1)\right] \le \tau_{ijk\theta \ttheta}, \quad \forall i \in I, j \in J, k \in K, \theta \in \Theta, \ttheta \in \hat{\Theta}, \label{queue2}\\
        & \sum_{m \in M_{ijk\ttheta}'} q_{ijkm\ttheta} \hat{x}_{ijkm \theta \ttheta} = \hat{\tau}_{ijk\theta \ttheta}, \quad \forall i \in I, j \in J, k \in K, \theta \in \Theta, \ttheta \in \hat{\Theta}, \label{queue3}\\
        & \mathbf{y} \in \mathcal{C}, \\
        & x_{ijkm\theta\ttheta},\hat{x}_{ijkm\theta\ttheta} \in \{0,1\}, \quad \forall i \in I, j \in J, k \in K, \theta \in \Theta,\ttheta \in \hat{\Theta}, m \in M{ijk\hat{\Theta}}, \\
        & y_{ijk\theta}, \hat{y}_{ijk \theta \ttheta} \in \{0,1\}, \quad \forall i \in I, j \in J, k \in K, \theta \in \Theta,\ttheta \in \hat{\Theta}, \\
        & z_{ikn} \in \{0,1\}, \quad \forall i \in I, k \in K, n \in {1,...,N_k}, \\
        & \tau_{ijk\theta \ttheta} \ge 0, \quad \forall i \in I, j \in J, k \in K, \theta \in \Theta, \ttheta \in \hat{\Theta}. 
        \end{align}
    \section{Rural Uncertainty Parameter}\label{app_rural_analysis}
    The rural dataset consists of 392,946 unique EMS responses, which we divided into training (196,473 responses), testing (98,236 responses), and simulation (98,237 responses) subsets. Following the same methodology applied to urban data, we determined a classification threshold of 0.689 to identify patients eligible for diversion from emergency departments.
    \begin{figure}[h! tbp]
        \centering
        \includegraphics[width=0.8\linewidth]{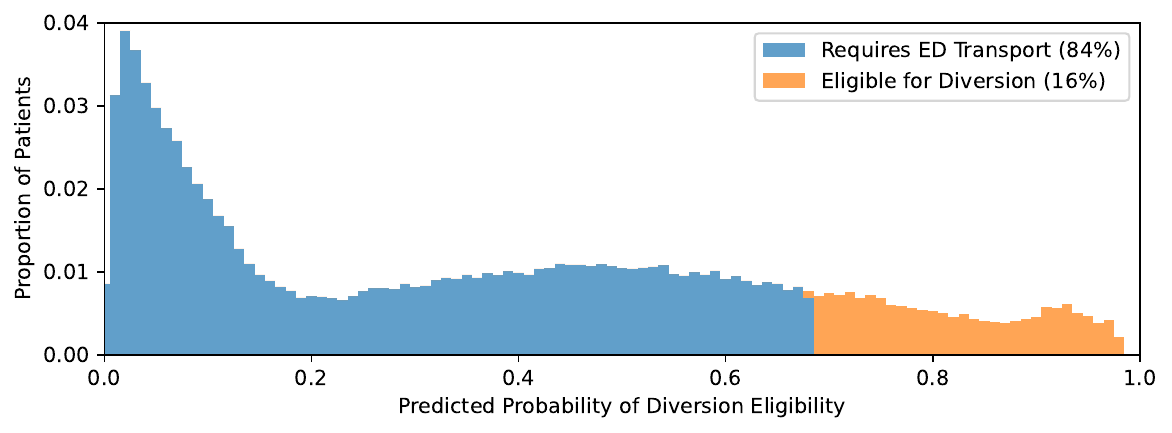}
        \caption{Distribution of predicted probabilities from the first-stage model used to estimate diversion eligibility in rural environments. A classification threshold of 0.689 was selected to reflect prior estimates suggesting that approximately 16\% of rural EMS patients may be safely diverted.}        \label{fig:rural_pred_prob}
    \end{figure}
    \begin{figure}[h! tbp]
        \centering
        \includegraphics[width=0.8\linewidth]{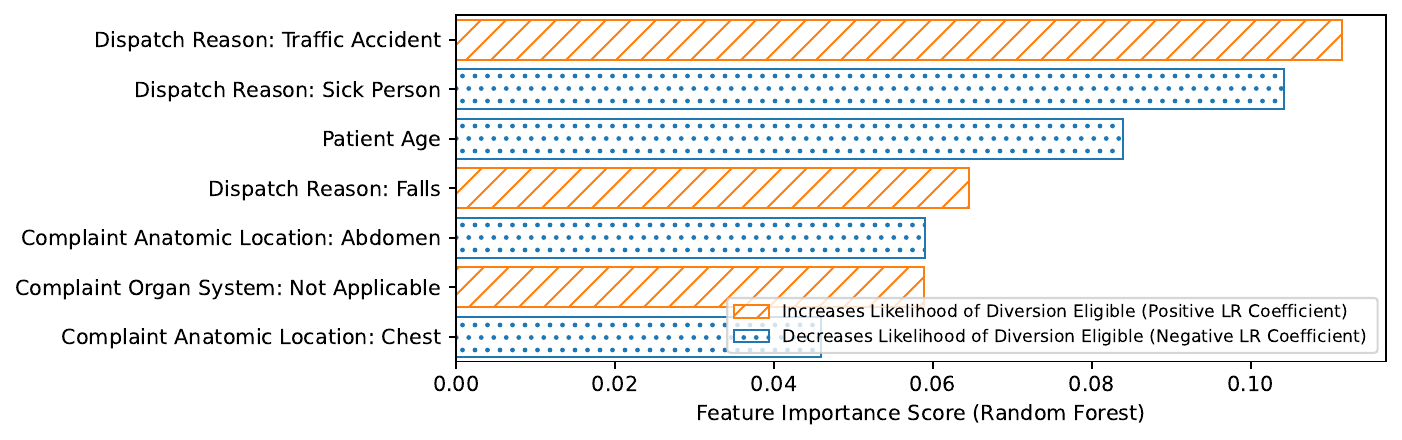} 
        \caption{Feature importance from a random forest model trained on phone-screening data for rural environments. Hatching denotes the direction of association with diversion eligibility, based on the sign of the corresponding coefficient in a logistic regression model.}    \label{fig:rural_features}
    \end{figure}\\
    Table \ref{tab:screening_accuracy_rural} summarizes the probability parameters $p_{\theta \ttheta}$ for rural environments used in our optimization experiments. Similar to the urban analysis, these parameters capture the uncertainty in dispatcher screening and inform our allocation model for rural PC-EMS implementation. 
    \begin{table}[h! tbp]
    \centering
    \renewcommand{\arraystretch}{.7} 
    \caption{Estimated probability of actual treatment needs ($\hat{\Theta}$) conditioned on screening predictions ($\Theta$). These probabilities define the parameter $p_{\theta \ttheta}$ in our optimization model.}
    \label{tab:screening_accuracy_rural}
    \begin{tabular}{|C{6cm}|C{3.1cm}||C{1.6cm}|C{1.6cm}|C{1.6cm}|}
    \hline
    \multirow{2}{*}{\shortstack[c]{\textbf{Believed Condition ($\Theta$)}}} & 
    \multirow{2}{*}{\raisebox{-.25\height}{\shortstack[c]{\textbf{\% of Patients}\\ \textbf{Screened}}}} & \multicolumn{3}{c|}{\textbf{Actual Condition ($\hat{\Theta}$)}} \\ 
    \cline{3-5}
    & & \textbf{ED} & \textbf{AD} & \textbf{TIP} \\
    \hline
    Likely Not Diversion Eligible (ED)  & 68.9\%  & 94.2\%  & 0.4\%  & 5.4\%  \\
    \hline
    Likely Diversion Eligible (AD/TIP)  & 31.1\%  & 61.2\%  & 2.6\%  & 36.2\%  \\
    \hline
    \end{tabular}
    \end{table}
    \section{Dispatching Strategies and Number of Diversion-Capable Units} \label{app_dispatch_sens}
    Figure \ref{fig:app_cont_comp} compares different dispatching strategies across varying proportions of diversion-capable units in the fleet. While the main paper presents results for specific configurations (10\%, 25\%, and 50\%), Figure \ref{fig:app_cont_comp} illustrates the continuous relationship between fleet composition and diversion rates for all three dispatching strategies across the four geographic environments.
    \begin{figure}[h! tbp]
        \centering
        \includegraphics[width=\linewidth]{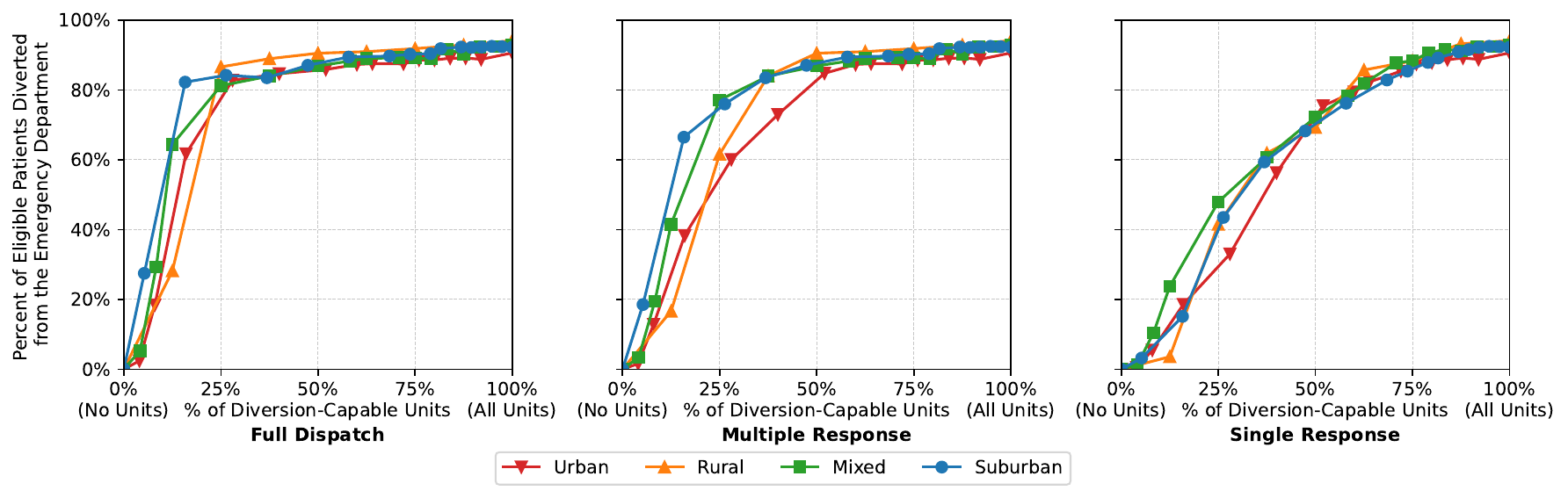}
        \caption{Impact of dispatching strategies on diversion rate.}
        \label{fig:app_cont_comp}
    \end{figure}
    \section{Screening Accuracy and Number of Diversion-Capable Units} \label{app_screening_sens}
    Figure \ref{fig:dis_cap} examines the interaction between patient screening accuracy and the proportion of diversion-capable units across different EMS environments. Section \ref{sec:value_screening} presented results for the 25\% fleet composition case, here we extend that analysis to include 10\% and 50\% diversion-capable fleet compositions to provide a more comprehensive understanding of these dynamics.
    \begin{figure}[ht!]
        \centering
        \subfloat[$\sim 10\%$ diversion-capable units]{%
            \includegraphics[width=0.8\textwidth]{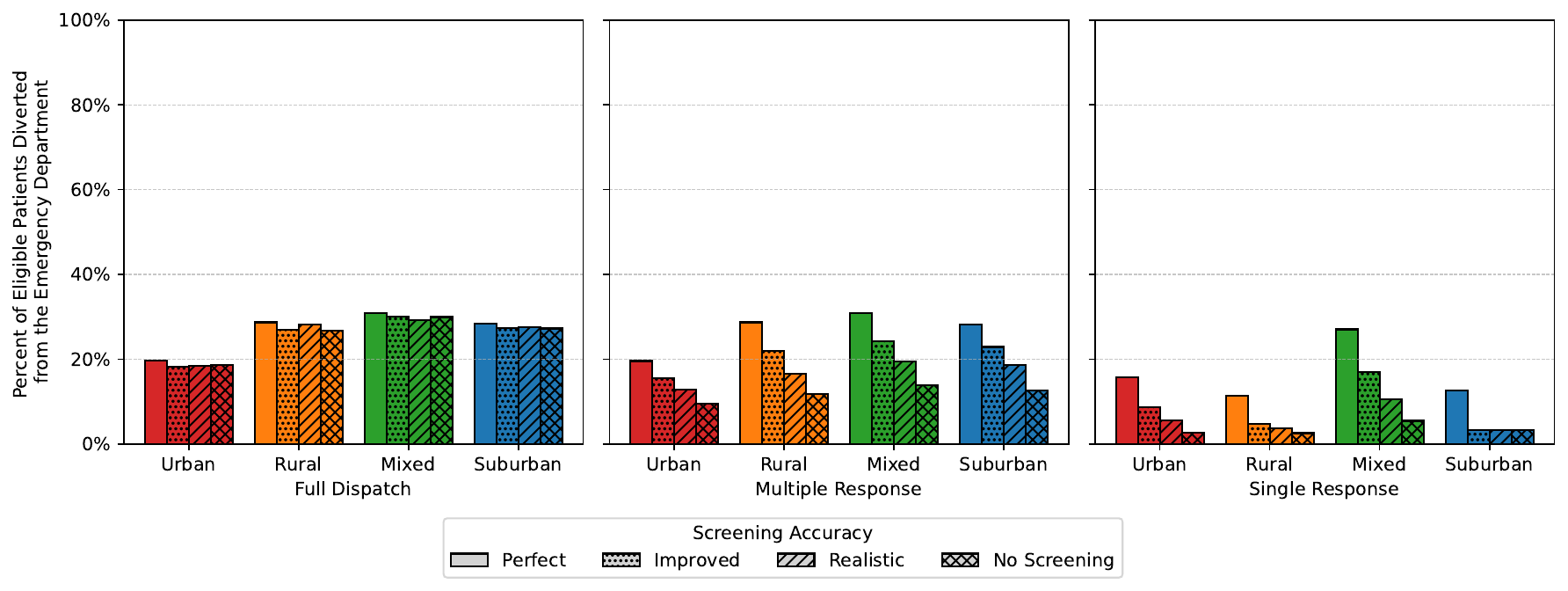}%
            \label{fig:dis-cap-10pct}%
        }
        
        \vspace{0.5em}
        
        \subfloat[$\sim 25\%$ diversion-capable units]{%
            \includegraphics[width=0.8\textwidth]{figure_pdfs/improved_visualization_25percent.pdf}%
            \label{fig:dis-cap-25pct}%
        }
        
        \vspace{0.5em}
        
        \subfloat[$\sim 50\%$ diversion-capable units]{%
            \includegraphics[width=0.8\textwidth]{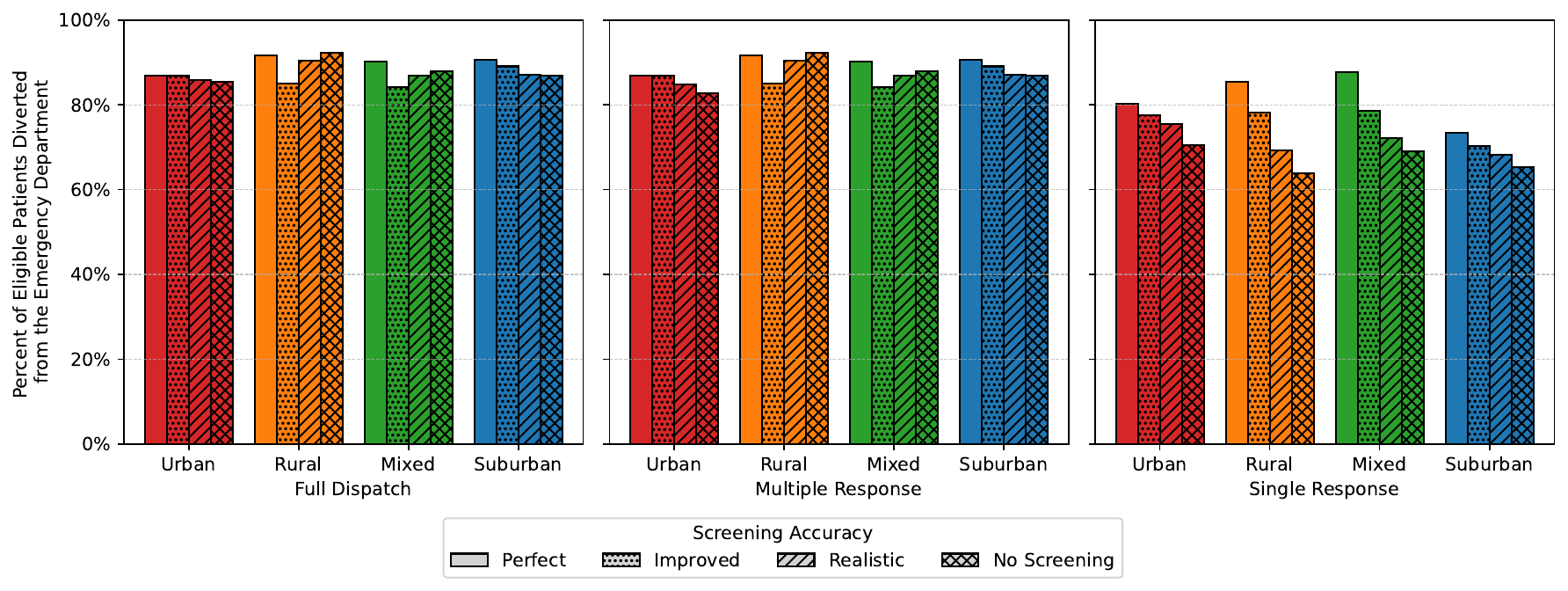}%
            \label{fig:dis-cap-50pct}%
        }
        \caption{Effect of patient screening accuracy on optimal diversion rates with varying diversion-capable fleet compositions across different geographic environments and dispatching strategies.}
        \label{fig:dis_cap}
    \end{figure}
 \end{APPENDICES}

\end{document}